\documentclass[a4paper,11pt,oneside]{article}
\usepackage[utf8]{inputenc}

\usepackage{amsmath,amssymb}
\usepackage{algorithm, algo_olm}
\usepackage{amsthm, mathtools}
\usepackage{amsfonts}
\usepackage{rotating}

\usepackage{tikz}
\usepackage{pgfplots}

\usepackage[english]{babel} 
\usepackage[utf8]{inputenc}
\usepackage{geometry}
\usepackage{multicol}
\usepackage{graphicx} 
\newtheorem{thm}{Theorem}[section]

\newtheorem{proper}[thm]{Property}

\newlength\figureheight
\newlength\figurewidth

\title{A Fast Summation Method for translation invariant kernels} 

\author{Fabien Casenave\thanks{IGN LAREG, Univ Paris Diderot, Sorbonne Paris Cit\'e, 5 rue Thomas Mann, 75205 Paris Cedex 13, France.}}

\begin{document}
\maketitle
\date{}

\begin{abstract}
We derive a Fast Multipole Method (FMM) where a low-rank approximation of the kernel is obtained using the Empirical Interpolation Method (EIM). 
Contrary to classical interpolation-based FMM, where the interpolation points and basis are fixed beforehand, 
the EIM is a nonlinear approximation method which constructs interpolation points and basis which are adapted to the kernel under consideration.
The basis functions are obtained using evaluations of the kernel itself.
We restrict ourselves to translation-invariant kernels, for which a modified version of the EIM approximation can be used in a multilevel FMM context;
we call the obtained algorithm Empirical Interpolation Fast Multipole Method (EIFMM).
An important feature of the EIFMM is a built-in error estimation of the interpolation error made by the low-rank approximation of the far-field behavior of the kernel:
the algorithm selects the optimal number of interpolation points required to ensure a given accuracy for the result, 
leading to important gains for inhomogeneous kernels.
\end{abstract}


\section{Introduction}
\label{intro}

We consider the following problem: compute approximations of many sums of the form
\begin{equation}
\label{eq:FMMsum}
f(\bar{x}_i)=\sum_{j=1}^N \sigma_j K(\bar{x}_i,\bar{y}_j), 1\leq i\leq N,
\end{equation}
for a set of potentials $\sigma_j$, and target points $\bar{x}_i\in\Omega$ and source points $\bar{y}_j\in\Omega$, where $\Omega\subset\mathbb{R}^3$,
with a controllable error $\epsilon$. 
The domain $\Omega$ and $N$ are fixed for all computations, but the potentials $\sigma_j$ and the points $\bar{x}_i$ and $\bar{y}_j$ can vary from one computation to another.
The many-query situation can come from a temporal N-body problem (in which case the set
of source points equals the set of observation points), or the iterative resolution of a dense linear system.

Various algorithms have been proposed to tackle such problems.
Some methods only using numerical evaluations of $K$ have been proposed, see~\cite{wavelets, lowrank, coulomb}.
Many methods, e.g. the mosaic skeleton method~\cite{tyrtyshnikov}, 
hierarchical matrices~\cite{hackbusch}, the Fast Multipole Method (FMM)~\cite{Greengard} (and its algebraic generalization H2-matrices~\cite{hackbusch3, hackbusch4}),
the panel clustering method~\cite{hackbusch2}, rely on local low-rank approximations of $K$:
\begin{equation}
\label{eq:locallowrank}
K(x,y)\approx\sum_{l,m=1}^d \Delta_{l,m}g_l(x)h_m(y), x\in I, y\in J,
\end{equation}
where $I$ and $J$ are subsets of $\Omega$. Such approximations can be analytic, like interpolation (for instance the Chebyshev interpolation as in
\cite{giebermann, hackbusch}), or multipole expansions.
FMM for general kernels have been developed leading to kernel-independent procedures, see~\cite{genFMM1, genFMM2}.
In~\cite{bbFMM}, another kernel-independent FMM algorithm, based on the Chebyshev interpolation to obtain the low-rank approximation~\eqref{eq:locallowrank} of $K$,
is developed and called Black-Box FMM.
Since the algorithm we present here is in the same family as the one developed in~\cite{bbFMM}, our numerical experiments contain comparisons
with this method.
For other kernel-independent interpolation based procedures, see~\cite{cheby0,fft}.
Adaptive Cross Approximation (ACA)~\cite{ACA, Hmat} directly constructs approximations from the evaluations of the kernel on the target and source points $K(\bar{x}_i,\bar{y}_j)$,
by selecting some columns and rows of the matrix $\left(K(\bar{x}_i,\bar{y}_j)\right)_{i,j}$.
Approximating a matrix by extracting some of its rows and columns has been investigated in~\cite{goreinov}.
With ACA methods,~\eqref{eq:locallowrank} is such that $g_l$ depends on $J$ and $h_m$ depends on $I$, which leads to a storage complexity of $dN\log(N)$.
With FMM and H2-matrices, $g_l$ and $h_m$ are independent of respectively $J$ and $I$, which allows additional factorization and the use of nested basis to lower the
storage complexity to $dN$. An improved ACA exploiting the nested basis property has been proposed in~\cite{bebendorf}.

In this work, a low-rank approximation of $K$ is obtained with the Empirical Interpolation Method (EIM), a procedure that represents a two-variable
function in a so-called affine dependent formula and that is used in the Reduced Basis community,~\cite{Barrault}. The EIM 
provides a non-linear and function-dependent interpolation: the basis functions and interpolation points are not fixed a priori but are computed during a greedy procedure,
which differs from the Chebyshev interpolations.
The procedure should capture any irregularity of the kernel better than a
fixed basis interpolation procedure, since the basis functions of the EIM are based on evaluations of the kernel itself.
Theoretical results on the convergence of EIM interpolation error in the general case
are limited, to the author knowledge, to an upper bound of $2^d-1$ of the Lebesgue constant, where $d$ is the number of interpolation points, see~\cite[Proposition 3.2]{Barrault}.
Various numerical studies illustrate the efficiency of the EIM~\cite{EIMeff3, EIMeff2, Maday, EIMeff1}, even though a theoretical analysis validating this efficiency is not available yet.
Depending on the kernel (and on the level in the tree for inhomogeneous kernels), the EIM selects an optimal number of terms that ensures a given
accuracy of the approximation. Such a variable-order interpolation procedure have been used in~\cite{borm} in the H2-matrix context with Lagrange polynomials.

A direct application of the EIM algorithm to the kernel $K$ to obtain a local low-rank approximation~\eqref{eq:locallowrank} is such that $g_l$ depends on $J$ and $h_m$ depends on $I$,
and therefore is not compatible with a multilevel FMM: a nested basis procedure cannot be derived.
To correct this, we make use of the translation-invariant hypothesis of $K$ and derive
a local-global low-rank EIM approximation of $K$, in the sense that one variable is restricted to a local small subdomain while the
other variable is taken in almost all the domain. 
Then, we derive recursive formulae to produce a multilevel FMM for this modified approximation: this is detailed in Section~\ref{sec:multilevel}.
The modified approximation uses four summations, instead of two summations for the Black-Box FMM.
However, we will see that the M2L step, the most expensive one, has the same complexity as the Black-Box FMM one with respect to the number on interpolation points.
Our algorithm requires a precomputation step, during which the greedy step of the EIM at each level is carried-out, at well as compression of some operators.
In our simulations, we did not notice any speedup on homogeneous kernels compared to Black-Box FMM.
For inhomogeneous kernels, the approximation problem gets easier as we get higher in the tree: corresponding EIM require much less interpolation points and speedups are measured.
In $10^6$ points test-cases, depending on the kernel, the EIFMM gets more interesting if one has to compute more than 2 to 18 queries of the summation.
In a $10^8$ points test-case, for the kernel $e^{-r^2}$, the precomputation step becomes negligible, and a single query is computed twice as fast for a nearly one order of magnitude better accuracy.

In Section~\ref{sec:synthEIM}, the EIM is recalled. In Section~\ref{EIFMM} is presented the new FMM algorithm, called the Empirical Interpolation Fast Multipole Method (EIFMM):
first, a four-summation formula based on the EIM is derived in Section~\ref{sec:suitable_approx}, that can be used in a multilevel FMM procedure, then
a monolevel version of the EIFMM is proposed in Section~\ref{sec:monolevel}, and a multilevel one in Section~\ref{sec:multilevel}.
In Section~\ref{sec:optimization}, the overall execution time is reduced using some precomputation and classical compression of the M2L step.
Finally, numerical experiments are presented in Section~\ref{sec:num}, and conclusions are drawn in Section~\ref{sec:conclusion}.

\section{The Empirical Interpolation Method}
\label{sec:synthEIM}

Consider a function $K(x,y)$ defined over $\mathcal{D}_x\times\mathcal{D}_y$ assumed to be real-valued for simplicity, where $\mathcal{D}_x, \mathcal{D}_y\subset\Omega$.
Fix an integer $d$.
The Empirical Interpolation Method provides a way to approximate this function in the following form:
\begin{equation}
\label{eq:approx_EIM}
K(x,y) \approx (I_d K)(x,y):=\sum_{m=1}^d \lambda_m(x)q_m(y),
\end{equation}
where $\lambda_m(x)$ is such that
\begin{equation}
\label{eq:onlineapb}
\sum_{m=1}^d B_{l,m}\lambda_m(x) = K(x,y_l), \qquad\forall 1\leq l\leq d.
\end{equation}
The functions $q_{m}(\cdot)$ and the matrix $B\in\mathbb{R}^{d\times d}$, which is
lower triangular with unity diagonal, are constructed
in Algorithm~\ref{algo0}, where $\delta_d={\rm Id}-I_{d}$ and $\|\cdot\|_{\mathcal{D}_y}$ is a norm
on $\mathcal{D}_y$, for instance the $L^\infty\left(\mathcal{D}_y\right)$- or the $L^2\left(\mathcal{D}_y\right)$-norm.
In practice, the argmax appearing in Algorithm~\ref{algo0} is searched over finite subsets of
$\mathcal{D}_x$ and $\mathcal{D}_y$, denoted respectively by $\mathcal{D}_{x, {\rm trial}}$ and $\mathcal{D}_{y, {\rm trial}}$ and called training sets.
Notice that Algorithm~\ref{algo0} also constructs the set of points $\{y_l\}_{1\leq l\leq d}$ in $\mathcal{D}_y$ used in~\eqref{eq:onlineapb},
and a set of points $\{x_l\}_{1\leq l\leq d}$ in $\mathcal{D}_x$. These points are interpolations points for the approximation of $K$ (see Property~\ref{interp1})
and are different from the target and source points $\bar{x}_i,\bar{y}_j$.

\begin{algorithm}[h]
	\caption{Greedy algorithm of the EIM}
	\label{algo0}
	\begin{algorithmic}[1]
        \STATE {Choose $d>1$}
	\STATE {Set $k:=1$}
        \STATE {Compute $\displaystyle x_1:=\underset{x\in\mathcal{D}_x}{\textnormal{argmax}}\|K(x,\cdot)\|_{\mathcal{D}_y}$}
        \STATE {Compute $\displaystyle y_1:=\underset{y\in\mathcal{D}_y}{\textnormal{argmax}}|K(x_1,y)|$}
        \hfill \COMMENT{First interpolation point}
	\STATE {Set $\displaystyle q_1(\cdot):=\frac{K(x_1,\cdot)}{K(x_1,y_1)}$}
       \hfill \COMMENT{First basis function}
        \STATE {Set $B_{1,1}:=1$}
       \hfill \COMMENT{Initialize matrix $B$}
        \WHILE {$k < d$} 
		\STATE Compute $\displaystyle x_{k+1}:=\underset{x\in\mathcal{D}_x}{\textnormal{argmax}}\|(\delta_k K)(x,\cdot)\|_{\mathcal{D}_y}$
                \STATE Compute $\displaystyle y_{k+1}:=\underset{y\in\mathcal{D}_y}{\textnormal{argmax}}|(\delta_k K)(x_{k+1},y)|$
                \hfill \COMMENT{$(k+1)$-th interpolation point}  
                \STATE Set $\displaystyle q_{k+1}(\cdot):=\frac{(\delta_k K)(x_{k+1},\cdot)}{(\delta_k K)(x_{k+1},y_{k+1})}$
                \hfill \COMMENT{$(k+1)$-th basis function}
                \STATE Set $\displaystyle B_{k+1,m}:=q_{m}(y_{k+1})$, for all $1\leq m\leq {k+1}$
                \hfill \COMMENT{Increment matrix $B$}
                \STATE $k\leftarrow k+1$
                \hfill \COMMENT{Increment the size of the decomposition}
	\ENDWHILE
\end{algorithmic}
\end{algorithm}

The following assumption is made:
\begin{itemize}
 \item[\bf(H)] The dimension of $\underset{x\in\mathcal{D}_x}{\rm Span}\left(K(x,\cdot),\right)$ is larger than $d$.
\end{itemize}
From [{\bf(H)}], the functions
$\{K(x_l,\cdot)\}_{1\leq l\leq d}$ are linearly independent (otherwise,
$(\delta_k K)(\mu_{k+1},x_{k+1})=0$ for some $k$ in Algorithm~\ref{algo0}).

It is shown in~\cite{casenave_ACM} that the approximation~\eqref{eq:approx_EIM} can be written in the form
\begin{equation}
\label{eq:approx2}
\begin{aligned}
\left(I_{d} K\right)(x,y) = \sum_{l,m=1}^{d} \Delta_{l,m} K(x,y_l)K(x_m,y),
\end{aligned}
\end{equation}
where $\Delta = B^{-t}\Gamma^{-1}$, and where
the matrix $\Gamma$ is is upper triangular and constructed recursively in the loop in $k$ of~Algorithm~\ref{algo0} in the following way:
\begin{itemize}
 \item $k=1$:
\begin{equation*}
\Gamma_{1,1}=K(x_{1}, y_{1}), 
\end{equation*}
 \item $k\rightarrow k+1$: 
\begin{equation*}
\begin{alignedat}{3}
\Gamma_{k+1,k+1}&=(\delta_{k} K)(x_{k+1}, y_{k+1}),&&\\ 
\Gamma_{m,k+1}&=0,&\qquad \forall 1\leq m\leq k,&\\
\Gamma_{k+1,m}&=\alpha_{m},&\qquad \forall 1\leq m\leq k,&\\ 
\end{alignedat}
\end{equation*}
where the vector $\alpha$ is such that $\sum_{m=1}^k B_{l,m}\alpha_{m}=K(x_{k+1}, y_{l})$, for all $1\leq~l\leq~k$.
\end{itemize}
The particular form~\eqref{eq:approx2} is the key to the use of the EIM approximation in a fast multipole context.

We recall the interpolation property of $I_d K$; see~\cite[Lemma~1]{Maday}:
\begin{proper}[Interpolation property]
\label{interp1}
For all $1\leq m\leq d$,
\begin{equation*}
\left\{
\begin{aligned}
(I_d K)(x,y_m) &= K(x,y_m), \quad \textnormal{for all } x\in\mathcal{D}_x,\\
(I_d K)(x_m,y) &= K(x_m,y), \quad \textnormal{for all } y\in\mathcal{D}_y.
\end{aligned}
\right.
\end{equation*}
\end{proper}

An important aspect is that the error made by the EIM approximation over the training set $\mathcal{D}_{x, {\rm trial}}\times\mathcal{D}_{y, {\rm trial}}$,
denoted $\epsilon_k$ and defined as $\epsilon_k:=\underset{y\in\mathcal{D}_y}{\textnormal{max}}|(\delta_k K)(x_{k+1},y)|$, is monitored during
the greedy procedure.
In practice, a stopping criterion is used on $\epsilon_k$ instead of on the number of interpolation points $d$, which was used to simplify the presentation of the EIM.
Hence, depending on the regularity of the function $K$, the EIM will automatically select the number of interpolation points
to certify the error over $\mathcal{D}_{x, {\rm trial}}\times\mathcal{D}_{y, {\rm trial}}$.

\section{The Empirical Interpolation Fast Multipole Method}
\label{EIFMM}

The FMM algorithm is based on a low-rank approximation of the far-field behavior of the kernel of the form~\eqref{eq:locallowrank}.
Such approximations enable some factorizations in the computation of key-quantities and the propagation of information through a tree structure, resulting
in lowering the complexity of the computation of~\eqref{eq:FMMsum}.
For a general presentation of the FMM algorithm, see the original article~\cite[Section 3]{Greengard}.

\subsection{A suitable approximation}
\label{sec:suitable_approx}

For the simplicity of the presentation, the algorithm is described in 1D. It can be readily extended in 2D and 3D.
Consider a kernel $K(x,y)$, $(x,y)\in\Omega\times \Omega$, where $\Omega=(-\frac{L}{2},\frac{L}{2})$ is the complete domain on which the FMM computation is carried out.
The kernel is supposed translation invariant, i.e. such that $K(x,y)=K(x-a,y-a)$, $\forall a\in\mathbb{R}$.
The domain $\Omega$ is partitioned in $2^\kappa$ intervals, $\kappa\in\mathbb{N}^*$, of length $2l_\kappa$, with $l_\kappa:=2^{-\kappa-1} L$.
Consider two intervals $I$ and $J$, well-separated in the sense that $\underset{(x,y)\in I\times J}{\min}|x-y|\geq 2l_\kappa$.
To anticipate the generalization to higher dimensions, $I$ and $J$ are called ``boxes'' in what follows.
We denote $c_I$ and $c_J$ the center of the boxes $I$ and $J$. 
Consider an EIM approximation constructed over $I\times J$:
\begin{equation}
\label{eq:approx3}
\begin{aligned}
K(x,y) \approx \sum_{l,m=1}^{d} \Delta^{I,J}_{l,m} K(x,y_l^J)K(x_m^I,y),\qquad x\in I,\quad y\in J,
\end{aligned}
\end{equation}
where $I$ and $J$ superscripts have been added to underline the box-dependent quantities (as deduced from Algorithm~\ref{algo0}).
In~\eqref{eq:approx3}, $K(x,y_l^J)$ depends on $J$, the box containing $y$ and $K(x_m^I,y)$ depends on $I$, the box containing $x$:
as stated in the introduction, this approximation is then not suitable for a multilevel FMM procedure.

\begin{figure}[h!]
\input{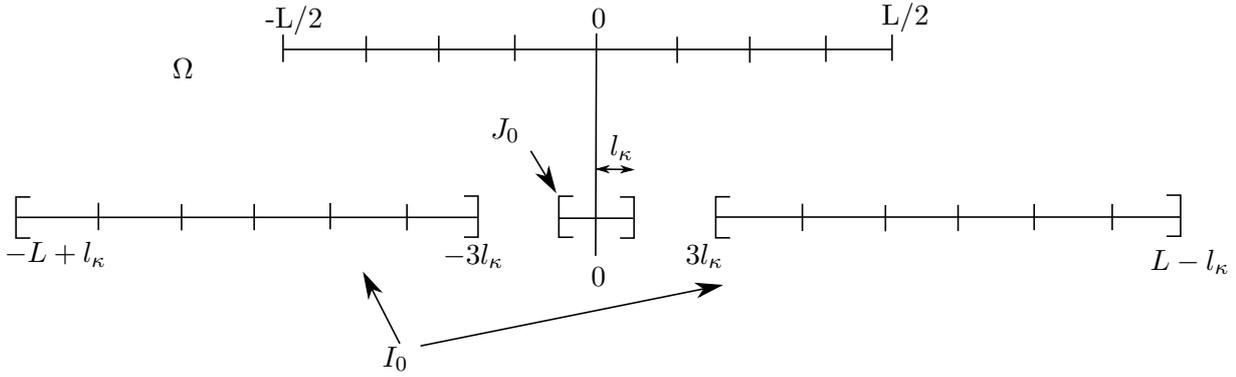}
 \caption{Representation of the partitioning of $\Omega$ in 1D for $\kappa=3$ and of the corresponding boxes $I_0$ and $J_0$.}
\label{fig:box_def}
\end{figure}

Define now $I_0:=\{x-c_{\hat{J}}, x\in \hat{I}\}$ for all box $\hat{J}$ in the partitioning of $\Omega$
and all box $\hat{I}$ well separated from $\hat{J}$, and define $J_0$ the box of length $2l_\kappa$ located at the center of $\Omega$
(notice that $J_0$ is not a box from the partitioning of $\Omega$).
In our setting, $J_0=(-l_\kappa, l_\kappa)$ and $I_0=(-L+l_\kappa, -3l_\kappa)\cup (3l_\kappa,L-l_\kappa)$, see Figure~\ref{fig:box_def}.
In other words, $J_0$ is a source box translated to the center of the domain, $I_0$ is the union of possible target boxes, translated in the same fashion.

\vspace{0.5cm}

Since $K(x,y) =K(x-c_J,y-c_J)$, with $x-c_J\in I_0$ and $y-c_J\in J_0$,
it is inferred that
\begin{equation}
\label{eq:approx4}
\begin{aligned}
K(x,y) \approx \sum_{l,m=1}^{d} \Delta^{I_0,J_0}_{l,m} K(x-c_J,y_l^{J_0})K(x_m^{I_0},y-c_J),
\end{aligned}
\end{equation}
where the EIM is carried out on $\mathcal{D}_x=I_0$ and $\mathcal{D}_y=J_0$.
An important aspect is that $I_0$ and $J_0$ are fixed domains, independent from the boxes $I$ and $J$ containing respectively $x$ and $y$.
Consider the first evaluation of the kernel in~\eqref{eq:approx4}. There holds $K(x-c_J,y_l^{J_0})=K(x-c_I,y_l^{J_0}+c_J-c_I)$, where
$x-c_I\in J_0$ and $y_l^{J_0}+c_J-c_I\in I_0$. It is then inferred that
\begin{equation}
\label{eq:approx5}
\begin{aligned}
K(x-c_J,y_l^{J_0}) \approx \sum_{l',m'=1}^{d} \Delta^{J_0,I_0}_{l',m'} K(x-c_I,y_{l'}^{I_0})K(x_{m'}^{J_0},y_l^{J_0}+c_J-c_I),\qquad 1\leq l\leq m,
\end{aligned}
\end{equation}
where a second EIM has been carried out on $\mathcal{D}_x=J_0$ and $\mathcal{D}_y=I_0$ (with the same number of terms $d$ for simplicity).
Notice that if the kernel is symmetric, $\Delta^{I_0,J_0}=(\Delta^{J_0,I_0})^t$, $x_m^{I_0}=y_m^{I_0}$ for all $1\leq~m\leq~d$, and
$y_m^{J_0}=x_m^{J_0}$ for all $1\leq~m\leq~d$.
Injecting~\eqref{eq:approx5} into~\eqref{eq:approx4}:
\begin{equation}
\label{eq:approxFMM}
\begin{aligned}
K(x,y) \approx \sum_{l'=1}^{d} K(x-c_I,y_{l'}^{I_0})\sum_{m'=1}^{d} \Delta^{J_0,I_0}_{l',m'} \sum_{l=1}^{d}K(x_{m'}^{J_0},y_l^{J_0}+c_J-c_I)\sum_{m=1}^{d}\Delta^{I_0,J_0}_{l,m}K(x_m^{I_0},y-c_J),
x\in I,y\in J.
\end{aligned}
\end{equation}
In~\eqref{eq:approxFMM}, the only term depending on the boxes $I$ and $J$ (besides $x$ and $y$ themselves) are $c_I$ and $c_J$, the center of the boxes $I$ and $J$.
Hence, in the low-rank approximation, the $x$-dependent function $K(x-c_I,y_{l'}^{I_0})$ no longer depends on $J$ and the $y$-dependent function 
$K(x_m^{I_0},y-c_J)$ no longer depends on $I$. As a consequence, the non-trivial manipulation leading to~\eqref{eq:approxFMM} are required for the use
of EIM interpolation in a multilevel FMM context.

Notice that in all the kernel evaluations in~\eqref{eq:approxFMM}, the first and second variables are always separated by a distance
of at least $2l_{\kappa}$.

\subsection{A monolevel multipole method}
\label{sec:monolevel}

Consider a set of $N$ target points $\{\bar{x}_i\}_{1\leq i\leq N}$, $N$ source points $\{\bar{y}_j\}_{1\leq j\leq N}$ in~$\Omega$,
and $N$ potentials $\{\sigma_j\}_{1\leq j\leq N}$, with $N\gg 1$.
We look for a low complexity approximation of
\begin{equation}
f(\bar{x}_i) = \sum_{j=1}^N \sigma_j K(\bar{x}_i, \bar{y}_j), \qquad 1\leq i\leq N.
\end{equation}

We say that two boxes $I$ and $J$ are well-separated if $\underset{x,y\in I\times J}{\min}|x-y|\geq 2l_{\kappa}$.
The set of boxes well-separated from $I$ is denoted $\mathcal{F}(I)$, its complement in $\Omega$ is the neighborhood of $I$ and is denoted $\mathcal{N}(I)$.
The FMM aims at compressing the ``far'' interactions in the following fashion:
\begin{equation}
f_{\mathcal{F}}(\bar{x}_i) = \sum_{J\in\mathcal{F}(I)}\sum_{\bar{y}_j\in J} \sigma_j K(\bar{x}_i, \bar{y}_j), \qquad \textnormal{ for all box }I
\textnormal{ and all } \bar{x}_i\in I.
\end{equation}

A monolevel multipole method can be directly derived from~\eqref{eq:approxFMM} to approximate $f_{\mathcal{F}}(\bar{x}_i)$,
see Algorithm~\ref{algo1}. The obtained approximation of $f_{\mathcal{F}}(\bar{x}_i)$ is denoted $f_{\mathcal{F}, {\rm FMM}}(\bar{x}_i)$. 
An approximation of $f(\bar{x}_i)$ is then given by
\begin{equation}
f(\bar{x}_i)\approx f_{\mathcal{F}, {\rm FMM}}(\bar{x}_i)+f_{\mathcal{N}}(\bar{x}_i),
\end{equation}
where the ``near`` interactions are defined by
\begin{equation}
f_{\mathcal{N}}(\bar{x}_i) = \sum_{\bar{y}_j\in \mathcal{N}(I)} \sigma_j K(\bar{x}_i, \bar{y}_j).
\end{equation}
\begin{algorithm}[h!]
	\caption{Monolevel EIFMM}
	\label{algo1}
	\begin{algorithmic}[1]
        \STATE {Compute $W_m^I:=\sum_{\bar{y}_j\in I}\sigma_j K(x_m^{I_0},\bar{y}_j-c_I)$, for all box $I$ and all $1\leq m\leq d$}
        \STATE {Compute $\hat{W}_l^I:=\sum_{m=1}^{d}\Delta^{I_0,J_0}_{l,m}W_m^I$, for all box $I$ and all $1\leq l\leq d$}
        \STATE {Compute $l_{m'}^I:=\sum_{J\in\mathcal{F}(I)}\sum_{l=1}^{d}K(x_{m'}^{J_0},y_l^{J_0}+c_J-c_I)\hat{W}_l^J$
, for all box $I$ and all $1\leq m'\leq d$}
        \STATE {Compute $\hat{l}_{l'}^I:=\sum_{m'=1}^{d} \Delta^{J_0,I_0}_{l',m'}l_{l'}^I$, for all box $I$ and all $1\leq m'\leq d$}
        \STATE {Compute $f_{\mathcal{F}, {\rm FMM}}(\bar{x}_i):=\sum_{l'=1}^{d} K(\bar{x}_i-c_I,y_{l'}^{I_0})\hat{l}_{l'}^I$, for all box $I$ and all $x_i\in I$}
\end{algorithmic}
\end{algorithm}

\begin{table}
\begin{center}
   \begin{tabular}{| c | c |}
     \hline
     Step number & Complexity \\ \hline
     1 & $dN$ \\ \hline
     2 & $d^2 (2^D)^\kappa$ \\ \hline
     3 & $d^2 (2^D)^{2\kappa}$ \\ \hline
     4 & $d^2 (2^D)^\kappa$ \\ \hline
     5 & $dN$ \\     \hline
     $f_{\mathcal{N}}$ & $N^2 (2^D)^{-\kappa}$\\
     \hline
   \end{tabular}
 \end{center}
\caption{Complexity of each step of Algorithm~\ref{algo1} and of the near interactions term.}
\label{tab1}
\end{table}

In Table~\ref{tab1} is given the complexity of each step of~Algorithm~\ref{algo1}.
For step 3 and the computation of the near interactions to be in the same order in $N$, one must impose $(2^D)^{2^\kappa}\propto N^2 (2^D)^{-\kappa}$, leading to
a overall complexity of the sum of order $N^{\frac{4}{3}}$.
As classically done in the literature, the derivation of a multilevel multipole method allows to lower the complexity in $N$ of the algorithm to linear.

\subsection{A multilevel multipole method}
\label{sec:multilevel}

The domain of computation $\Omega=(-\frac{L}{2},\frac{L}{2})$ is now organized as a binary tree with $\kappa+1$ levels.
Level $0$ is the root of the tree: there is one box equal to $\Omega$, level $1$ contains the partitioning of
$\Omega$ into two boxes of length $\frac{L}{2}$, and the following levels are obtained by partitioning each box of the previous
level in the same fashion, up to level $\kappa$.

A superscript $k$ is now added to the quantities previously defined, to indicate that they depend on the level $k$.
Two EIM's are to be considered at each level, in the same fashion as in Section~\ref{sec:suitable_approx}, where now the domains depend on the level $k$ in the tree:
$I_0^k=(-L+l_k, -3l_k)\cup (3l_k, L-l_k)$, and $J_0^k=(-l_k, l_k)$. 

For symmetric kernels, this amount of precomputation can be lowered to one EIM by level. 
The key to a multilevel method stands in the ability to derive recursion formulae to propagate information
from the level $\kappa$ to the level $0$ in the upward pass, and from the level $0$ to the level $\kappa$ in the downward pass.

Black-Box FMM, at least in its original version, is proposed with a constant-order of interpolation, which
enables to propagate the information through the levels of the tree without introducing additional errors. Considering a variable-order Chebyshev interpolation is possible,
but one must carefully analyse the introduced error, as done in~\cite{borm} in the context of H2-matrices.
In our method, the domains of approximation depicted in Figure~\ref{fig:box_def} naturally enable a variable-order interpolation depending on the level of the tree,
for which the introduced error is the same as the EIM approximation error in~\eqref{eq:approxFMM}. As a consequence, the variable-order aspect of the
approximation is directly obtained at the required accuracy, and the length of the approximation is automatically optimized at each level of the tree, which leads to important
gains for inhomogeneous kernels.

\subsubsection{Recursion for the upward pass: multipole to multipole (M2M)}

The goal is to compute 
\begin{equation}
\label{eq:upward_pass}
{W}_m^{I^k}:=\sum_{\bar{y}_j\in I^k}\sigma_j K(x_m^{I_0^k},\bar{y}_j-c_{I^k}),
\end{equation}
for all box $I^k$ at all level $k$ with $0\leq k\leq\kappa$ and all $1\leq m\leq d^k$ by recursion. The initialization of the recursion consists in
computing ${W}_m^{I^\kappa}$ at level $\kappa$ using~\eqref{eq:upward_pass}. Then, suppose ${W}_m^{I^k}$ is available at level $k$, there holds
\begin{equation}
\label{eq:upward}
{W}_m^{I^{k-1}}=\sum_{\bar{y}_j\in I^{k-1}}\sigma_j K(x_m^{I_0^{k-1}},\bar{y}_j-c_{I^{k-1}})=\sum_{J^k\in\mathcal{C}(I^{k-1})}\sum_{\bar{y}_j\in J^k}\sigma_j 
K(x_m^{I_0^{k-1}}+c_{I^{k-1}}-c_{J^k},\bar{y}_j-c_{J^k}),
\end{equation}
where $\mathcal{C}(I^{k-1})$ denotes the children boxes of $I^{k-1}$.
To obtain a recursion on ${W}_m^{I^{k-1}}$, i.e. to compute ${W}_m^{I^{k-1}}$ using some values of ${W}_m^{J^{k}}$, for $J^k$ being boxes at the level $k$, we will
use the EIM approximations already computed at level $k$ to derive~\eqref{eq:approxFMM}. One must check that the corresponding EIM approximations are evaluated on the right
domains, namely $I_0^k$ and $J_0^k$.

For all $\bar{y}_j\in J^k$, $-l_k\leq \bar{y}_j-c_{J^k}\leq l_k$, i.e. $\bar{y}_j-c_{J^k}\in J_0^k$.
Then, for all $J^k\in\mathcal{C}(I^{k-1})$, $-l_k\leq c_{I^{k-1}}-c_{J^k}\leq l_k$ (see for instance Figure~\ref{fig:closestIJ}).
Moreover, since $x_m^{I_0^{k-1}}\in I_0^{k-1}=(-L+l_{k-1}, -3l_{k-1})\cup (3l_{k-1}, L-l_{k-1})$, it is inferred that
\begin{equation*}
-L+l_k\leq x_m^{I_0^{k-1}}+c_{I^{k-1}}-c_{J^k}\leq -5l_k\leq -3l_k\textnormal{ or }3l_k\leq 5l_k\leq x_m^{I_0^{k-1}}+c_{I^{k-1}}-c_{J^k}\leq L-l_k,
\end{equation*}
i.e. $x_m^{I_0^{k-1}}+c_{I^{k-1}}-c_{J^k}\in I_0^k$.
Therefore, the EIM from level $k$ can then be used to write the following approximation:
\begin{equation}
\label{eq:EIM_upward}
K(x_m^{I_0^{k-1}}+c_{I^{k-1}}-c_{J^k},\bar{y}_j-c_{J^k})\approx \sum_{p,q=1}^{d^k}\Delta_{p,q}^{I_0^k, J_0^k}
K(x_m^{I_0^{k-1}}+c_{I^{k-1}}-c_{J^k},y_p^{J_0^k})K(x_q^{I_0^{k}},\bar{y}_j-c_{J^k}),
\end{equation}
where the approximation error is of the same order as in the~\eqref{eq:approxFMM}.

Then, injecting~\eqref{eq:EIM_upward} into~\eqref{eq:upward}, it is inferred that
\begin{equation}
\label{eq:rec_W}
\begin{aligned}
{W}_m^{I^{k-1}}&\approx\sum_{p,q=1}^{d^k}\Delta_{p,q}^{I_0^k, J_0^k}\sum_{J^k\in\mathcal{C}(I^{k-1})}K(x_m^{I_0^{k-1}}+c_{I^{k-1}}-c_{J^k},y_p^{J_0^k})\sum_{\bar{y}_j\in J^k}\sigma_j K(x_q^{I_0^{k}},\bar{y}_j-c_{J^k})\\
&\approx\sum_{p,q=1}^{d^k}\Delta_{p,q}^{I_0^k, J_0^k}\sum_{J^k\in\mathcal{C}(I^{k-1})}K(x_m^{I_0^{k-1}}+c_{I^{k-1}}-c_{J^k},y_p^{J_0^k})W_q^{J^k},
\end{aligned}
\end{equation}
which provides a recursive formula for the upward pass.
In the same fashion as for the monolevel version, we then compute
\begin{equation}
\label{eq:hatw}
\hat{W}_l^{I^k}:=\sum_{m=1}^{d^k}\Delta^{I_0^k,J_0^k}_{l,m}W_m^{I^k},
\end{equation}
for all box $I^k$ at all level k with $0\leq k\leq \kappa$ and all $1\leq l\leq d^k$.

\subsubsection{Transfer pass: multipole to local (M2L)}
\begin{figure}[h!]
\input{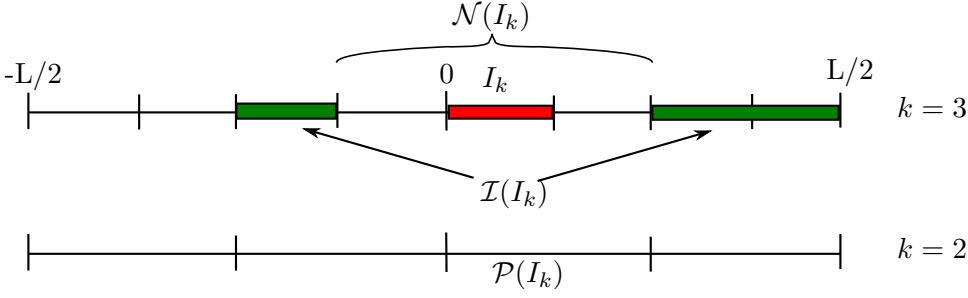}
 \caption{Representation of the neighborhood and the interaction list of a box $I^k$ and level $k$.}
\label{fig:defIk}
\end{figure}

We denote $\mathcal{I}(I^k)$, the set of boxes that are the children of the boxes in the neighborhood of the parent
of $I^k$ but are well-separated from $I^k$ at level $k$, and call it the interaction list of $I^k$, see Figure~\ref{fig:defIk}.

\vspace{0.5cm}

The M2L step consists in computing
\begin{equation}
\label{eq:M2Lstep}
{g}_{m'}^{I^k}:=\sum_{J^k\in\mathcal{I}(I^k)}\sum_{l=1}^{d^k}K(x_{m'}^{J_0^k},y_l^{J_0^k}+c_{J^k}-c_{I^k})\hat{W}_l^{J^k},
\end{equation}
for all box $I^k$ at all level k with $0\leq k\leq \kappa$ and all $1\leq m'\leq d^k$.

\subsubsection{Recursion for the downward pass: local to local (L2L)}

The goal is to compute 
\begin{equation}
\label{eq:downward_pass}
{l}_{m'}^{I^k}:=\sum_{J^k\in\mathcal{F}(I^k)}\sum_{l=1}^{d^k}K(x_{m'}^{J_0^k},y_l^{J_0^k}+c_{J^k}-c_{I^k})\hat{W}_l^{J^k},
\end{equation}
for all box $I^k$ at all level $k$ with $0\leq k\leq\kappa$ and all $1\leq m'\leq d^k$ by recursion.
The initialization of the recursion consists in
computing ${l}_{m'}^{I^2}$ at level $2$ using~\eqref{eq:downward_pass}. 
Then, suppose ${l}_{m'}^{I^k}$ is available, there holds
\begin{equation}
\label{eq:downward}
\begin{aligned}
{l}_{m'}^{I^{k+1}} &= {g}_{m'}^{I^{k+1}}+\sum_{J^{k}\in\mathcal{F}(\mathcal{P}(I^{k+1}))}\sum_{J^{k+1}\in\mathcal{C}(J^{k})}\sum_{l=1}^{d^{k+1}}K(x_{m'}^{J_0^{k+1}},y_l^{J_0^{k+1}}+c_{J^{k+1}}-c_{I^{k+1}})\hat{W}_l^{J^{k+1}}\\
&= {g}_{m'}^{I^{k+1}}+\sum_{J^{k}\in\mathcal{F}(\mathcal{P}(I^{k+1}))}\sum_{J^{k+1}\in\mathcal{C}(J^{k})}\sum_{l=1}^{d^{k+1}}K(x_{m'}^{J_0^{k+1}}+c_{I^{k+1}}-c_{\mathcal{P}(I^{k+1})},y_l^{J_0^{k+1}}+c_{J^{k+1}}-c_{\mathcal{P}(I^{k+1})})\hat{W}_l^{J^{k+1}},
\end{aligned}
\end{equation}
where $\mathcal{P}(I^{k+1})$ denotes the parent of $I^{k+1}$. 

\begin{figure}[h!]
\input{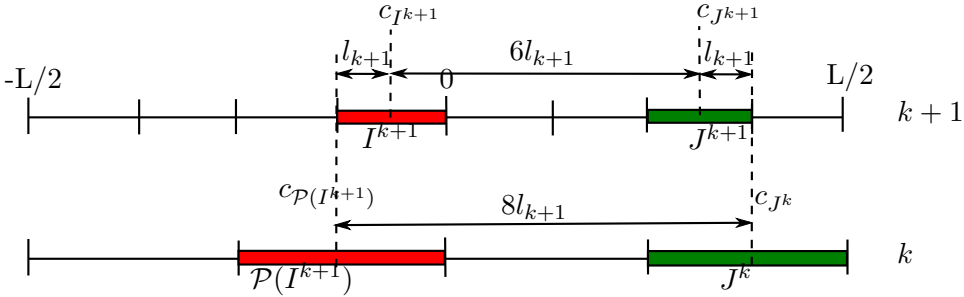}
 \caption{Possible values for $c_{I^{k+1}}-c_{\mathcal{P}(I^{k+1})}$ and $c_{J^{k+1}}-c_{J^{k}}$, and smallest possible values for $c_{J^{k+1}}-c_{\mathcal{P}(I^{k+1})}$
and $c_{\mathcal{P}(I^{k+1})}-c_{J^{k}}$, with $J^k$ well-separated from $\mathcal{P}({I^{k+1}})$ at level $k$, and $J^{k+1}$ child of $J^k$.}
\label{fig:closestIJ}
\end{figure}

\begin{figure}[h!]
\input{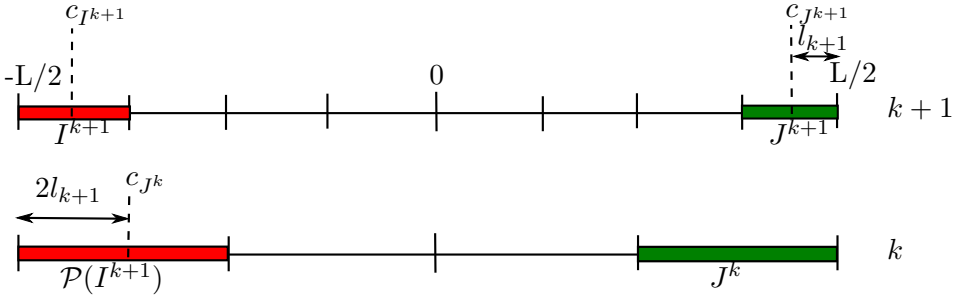}
 \caption{Largest possible value for $c_{J^{k+1}}-c_{\mathcal{P}(I^{k+1})}$ and $c_{\mathcal{P}(I^{k+1})}-c_{J^{k}}$, with $J^k$ well-separated from $\mathcal{P}({I^{k+1}})$ at level $k$, and $J^{k+1}$ child of $J^k$.}
\label{fig:farestIJ}
\end{figure}

It is seen on Figures~\ref{fig:closestIJ} and~\ref{fig:farestIJ} that
\begin{equation*}
-l_{k+1}\leq c_{I^{k+1}}-c_{\mathcal{P}(I^{k+1})}\leq l_{k+1},
\end{equation*}
and that
\begin{equation*}
-L+3l_{k+1}\leq c_{J^{k+1}}-c_{\mathcal{P}(I^{k+1})}\leq -7l_{k+1}\textnormal{ or }7l_{k+1}\leq c_{J^{k+1}}-c_{\mathcal{P}(I^{k+1})}\leq L-3l_{k+1}.
\end{equation*}

Since $x_{m'}^{J_0^{k+1}}, y_{l}^{J_0^{k+1}}\in J_0^{k+1}$, it
is then readily verified that $x_{m'}^{J_0^{k+1}}+c_{I^{k+1}}-c_{\mathcal{P}(I^{k+1})}\in J_0^{k}$ and
\begin{equation*}
-L+2l_{k+1}=-L+l_{k}\leq y_l^{J_0^{k+1}}+c_{J^{k+1}}-c_{\mathcal{P}(I^{k+1})}\leq -6l_{k+1}=-3l_k\textnormal{ or }3l_k\leq y_l^{J_0^{k+1}}+c_{J^{k+1}}-c_{\mathcal{P}(I^{k+1})}\leq L-l_{k},
\end{equation*}
i.e. $y_l^{J_0^{k+1}}+c_{J^{k+1}}-c_{\mathcal{P}(I^{k+1})}\in I_0^k$.
Therefore, the EIM from level $k$ can then be used to write the following approximation
\begin{equation}
\begin{aligned}
\label{eq:EIM_downward}
&\qquad\qquad K(x_{m'}^{J_0^{k+1}}+c_{I^{k+1}}-c_{\mathcal{P}(I^{k+1})},y_l^{J_0^{k+1}}+c_{J^{k+1}}-c_{\mathcal{P}(I^{k+1})})\approx\\
&\sum_{p',q'=1}^{d^k}\Delta_{p',q'}^{J_0^k, I_0^k} K(x_{m'}^{J_0^{k+1}}+c_{I^{k+1}}-c_{\mathcal{P}(I^{k+1})},y_{p'}^{I_0^k}) K(x_{q'}^{J_0^k},y_l^{J_0^{k+1}}+c_{J^{k+1}}-c_{\mathcal{P}(I^{k+1})}).
\end{aligned}
\end{equation}
Then, injecting~\eqref{eq:EIM_downward} into~\eqref{eq:downward}, it is inferred that
\begin{equation}
\begin{aligned}
\label{eq:relation_down}
{l}_{m'}^{I^{k+1}}\approx {g}_{m'}^{I^{k+1}}+&
\sum_{J^{k}\in\mathcal{F}(\mathcal{P}(I^{k+1}))}\sum_{J^{k+1}\in\mathcal{C}(J^{k})}
\sum_{p',q'=1}^{d^k}\Delta_{p',q'}^{J_0^k, I_0^k}\sum_{l=1}^{d^{k+1}}
K(x_{m'}^{J_0^{k+1}}+c_{I^{k+1}}-c_{\mathcal{P}(I^{k+1})},y_{p'}^{I_0^k})\\
& K(x_{q'}^{J_0^k},y_l^{J_0^{k+1}}+c_{J^{k+1}}-c_{\mathcal{P}(I^{k+1})})\hat{W}_l^{J^{k+1}}.
\end{aligned}
\end{equation}
Consider $K(x_{q'}^{J_0^k},y_l^{J_0^{k+1}}+c_{J^{k+1}}-c_{\mathcal{P}(I^{k+1})})=
K(x_{q'}^{J_0^k}+c_{\mathcal{P}(I^{k+1})}-c_{J^{k}},y_l^{J_0^{k+1}}+c_{J^{k+1}}-c_{J^{k}})$.
Since $J^{k}\in\mathcal{F}(\mathcal{P}(I^{k+1}))$ and $J^{k}$ is the parent of $J^{k+1}$, it is seen on Figures~\ref{fig:closestIJ} and~\ref{fig:farestIJ} that
\begin{equation*}
-l_{k+1}\leq c_{J^{k+1}}-c_{J^k}\leq l_{k+1},
\end{equation*}
and that
\begin{equation*}
-L+4l_{k+1}\leq c_{\mathcal{P}(I^{k+1})}-c_{J^{k}}\leq -8l_{k+1}\textnormal{ or }8l_{k+1}\leq c_{\mathcal{P}(I^{k+1})}-c_{J^{k}}\leq L-4l_{k+1}.
\end{equation*}
Since $x_{q'}^{J_0^k}, y_l^{J_0^{k+1}}\in J_0^{k+1}$, $J^{k}\in\mathcal{F}(\mathcal{P}(I^{k+1}))$, and $J^{k}$ is
the parent of $J^{k+1}$, it is readily verified that
$y_l^{J_0^{k+1}}+c_{J^{k+1}}-c_{J^{k}}\in J_0^{k}$ and
\begin{equation*}
-L+l_{k}\leq -L+3l_{k+1}\leq x_{q'}^{J_0^k}+c_{\mathcal{P}(I^{k+1})}-c_{J^{k}}\leq -7l_{k+1}\leq -3l_k\textnormal{ or }3l_k\leq x_{q'}^{J_0^k}+c_{\mathcal{P}(I^{k+1})}-c_{J^{k}}\leq L-l_{k},
\end{equation*}
i.e. $x_{q'}^{J_0^k}+c_{\mathcal{P}(I^{k+1})}-c_{J^{k}}\in I_0^k$.
Therefore, the EIM from level $k$ can then be used again to write the following approximation
\begin{equation}
\begin{aligned}
\label{eq:EIM_downward2}
&\qquad\qquad K(x_{q'}^{J_0^k}+c_{\mathcal{P}(I^{k+1})}-c_{J^{k}},y_l^{J_0^{k+1}}+c_{J^{k+1}}-c_{J^{k}})\approx\\
&\sum_{u,v=1}^{d^k}\Delta_{u,v}^{I_0^k, J_0^k} K(x_{q'}^{J_0^k}+c_{\mathcal{P}(I^{k+1})}-c_{J^{k}},y_{u}^{J_0^k}) K(x_{v}^{I_0^k},y_l^{J_0^{k+1}}+c_{J^{k+1}}-c_{J^{k}}).
\end{aligned}
\end{equation}
Injection~\eqref{eq:EIM_downward2} into~\eqref{eq:relation_down}, making use of the formula~\eqref{eq:hatw} to express
$\hat{W}_l^{J^{k+1}}$, it is inferred that
\begin{equation}
\begin{aligned}
\label{eq:relation_down2}
&{l}_{m'}^{I^{k+1}}\approx {g}_{m'}^{I^{k+1}}+
\sum_{J^{k}\in\mathcal{F}(\mathcal{P}(I^{k+1}))}\sum_{J^{k+1}\in\mathcal{C}(J^{k})}
\sum_{p',q'=1}^{d^k}\Delta_{p',q'}^{J_0^k, I_0^k}\sum_{l=1}^{d^{k+1}}
K(x_{m'}^{J_0^{k+1}}+c_{I^{k+1}}-c_{\mathcal{P}(I^{k+1})},y_{p'}^{I_0^k})\\
& \sum_{u,v=1}^{d^k}\Delta_{u,v}^{I_0^k, J_0^k} K(x_{q'}^{J_0^k}+c_{\mathcal{P}(I^{k+1})}-c_{J^{k}},y_{u}^{J_0^k}) 
K(x_{v}^{I_0^k},y_l^{J_0^{k+1}}+c_{J^{k+1}}-c_{J^{k}})
\sum_{m=1}^{d^{k+1}}\Delta^{I_0^{k+1},J_0^{k+1}}_{l,m}W_m^{J^{k+1}}.
\end{aligned}
\end{equation}
Making use of $K(x_{v}^{I_0^k},y_l^{J_0^{k+1}}+c_{J^{k+1}}-c_{J^{k}})=K(x_{v}^{I_0^k}+c_{J^{k}}-c_{J^{k+1}},y_l^{J_0^{k+1}})$
and $K(x_{q'}^{J_0^k}+c_{\mathcal{P}(I^{k+1})}-c_{J^{k}},y_{u}^{J_0^k})=K(x_{q'}^{J_0^k},y_{u}^{J_0^k}+c_{J^{k}}-c_{\mathcal{P}(I^{k+1})})$,
and reorganizing some terms, there holds
{\footnotesize
\begin{equation}
\begin{aligned}
\label{eq:relation_down3}
&{l}_{m'}^{I^{k+1}}\approx {g}_{m'}^{I^{k+1}}+\sum_{p',q'=1}^{d^k}\Delta_{p',q'}^{J_0^k, I_0^k}
K(x_{m'}^{J_0^{k+1}}+c_{I^{k+1}}-c_{\mathcal{P}(I^{k+1})},y_{p'}^{I_0^k})\\
&  \underbrace{\sum_{J^{k}\in\mathcal{F}(\mathcal{P}(I^{k+1}))}\sum_{u=1}^{d^k}K(x_{q'}^{J_0^k},y_{u}^{J_0^k}+c_{J^{k}}-c_{\mathcal{P}(I^{k+1})})
\underbrace{\sum_{v=1}^{d^k}\Delta_{u,v}^{I_0^k, J_0^k}
\underbrace{\sum_{l,m=1}^{d^{k+1}}\Delta^{I_0^{k+1},J_0^{k+1}}_{l,m}\sum_{J^{k+1}\in\mathcal{C}(J^{k})}K(x_{v}^{I_0^k}+c_{J^{k}}
-c_{J^{k+1}},y_l^{J_0^{k+1}})W_m^{J^{k+1}}}_{=W_v^{J^k}\textnormal{ using~\eqref{eq:rec_W}}}}_{=\hat{W}_u^{J^k}\textnormal
{ using~\eqref{eq:hatw}}}}_{=l_{q'}^{\mathcal{P}(I^{k+1})}\textnormal{ using~\eqref{eq:downward_pass}}},
\end{aligned}
\end{equation}}
\!\!\!\! which provides a recursive formula for the upward pass. In the same fashion as for the monolevel version, we then compute
\begin{equation}
\label{eq:hatl}
\hat{l}_{l'}^{I^\kappa}:=\sum_{m'=1}^{d^\kappa} \Delta^{J_0^\kappa,I_0^\kappa}_{l',m'}l_{l'}^{I^\kappa},
\end{equation}
for all box $I^\kappa$ at level $\kappa$ and all $1\leq l'\leq d^\kappa$.

\subsubsection{Algorithm and complexity}
 
Using the recursion formulae derived above, the multilevel procedure is detailed in Algorithm~\ref{algo2}.
\begin{algorithm}[h!]
	\caption{Multilevel EIFMM}
	\label{algo2}
	\begin{algorithmic}[1]
        \STATE {Compute ${W}_m^{I^\kappa}=\sum_{\bar{y}_j\in I^\kappa}\sigma_j K(x_m^{I_0^\kappa},\bar{y}_j-c_{I^\kappa})$, for all box $I^\kappa$ at level $\kappa$ and all $1\leq m\leq d^\kappa$}
        \STATE {Compute ${W}_m^{I^{k}}=\sum_{p,q=1}^{d^{k+1}}\Delta_{p,q}^{I_0^{k+1}, J_0^{k+1}}\sum_{J^{k+1}\in\mathcal{C}(I^{k})}K(x_m^{I_0^{k}}+c_{I^{k}}-c_{J^{k+1}},y_p^{J_0^{k+1}})W_q^{J^{k+1}}$,
for all box $I^k$ at all levels $\kappa-1\geq k\geq 0$ and all $1\leq m\leq d^{k+1}$}
         \STATE {Compute $\hat{W}_l^{I^k}=\sum_{m=1}^{d^k}\Delta^{I_0^k,J_0^k}_{l,m}W_m^{I^k}$,
for all box $I^k$ at all level k with $0\leq k\leq \kappa$ and all $1\leq l\leq d^k$}
        \STATE {Compute ${g}_{m'}^{I^k}=\sum_{J^k\in\mathcal{I}(I^k)}\sum_{l=1}^{d^k}K(x_{m'}^{J_0^k},y_l^{J_0^k}+c_{J^k}-c_{I^k})\hat{W}_l^{J^k}$,
for all box $I^k$ at all level k with $0\leq k\leq \kappa$ and all $1\leq m'\leq d^k$.}
        \STATE {Let ${l}_{m'}^{I^{0}}= {g}_{m'}^{I^{0}}$ and compute ${l}_{m'}^{I^{k}}= {g}_{m'}^{I^{k}}+\sum_{p',q'=1}^{d^{k-1}}\Delta_{p',q'}^{J_0^{k-1}, I_0^{k-1}}
K(x_{m'}^{J_0^{k}}+c_{I^{k}}-c_{\mathcal{P}(I^{k})},y_{p'}^{I_0^{k-1}})l_{q'}^{\mathcal{P}(I^{k})}$ for all box $I^k$ at all levels $1\leq k\leq \kappa$ and all $1\leq m'\leq d^{k-1}$}
        \STATE {Compute $\hat{l}_{l'}^{I^\kappa}:=\sum_{m'=1}^{d^\kappa} \Delta^{J_0^\kappa,I_0^\kappa}_{l',m'}l_{l'}^{I^\kappa}$ for all box $I^\kappa$ at level $\kappa$ and all $1\leq l'\leq d^\kappa$}
        \STATE {Compute $f_{\mathcal{F}, {\rm FMM}}(\bar{x}_i)=\sum_{l'=1}^{d^\kappa} K(\bar{x}_i-c_{I^\kappa},y_{l'}^{I_0^\kappa})\hat{l}_{l'}^{I^\kappa}$, for all box $I^\kappa$ and all $x_i\in I^\kappa$ at level $\kappa$}
\end{algorithmic}
\end{algorithm}

The complexity of the multilevel procedure is given in Table~\ref{tab22}.
For step 4 and the computation of the near interactions to be in the same order in $N$, one must impose $(2^D)^{\kappa}\propto N^2 (2^D)^{-\kappa}$, leading to
a overall complexity of the sum of order $N$.
The choice that minimizes the overall cost is $\kappa\propto \log(N)$ and a constant number of points per box at the level $\kappa$. 
As for any multilevel FMM algorithm, the complexity of step 3 in Algorithm~\ref{algo1} has been reduced by replacing the sum over $J\in\mathcal{L}(I)$ (which
is of complexity linear in $N$) by a sum over $J\in\mathcal{I}(I)$ (which is of complexity independent in $N$),
making use of the fact that the set of boxes well-separated from $I^k$ at level $k$ contains the set of boxes
well-separated from the parent of $I^{k}$ at level $k-1$.
By precomputing some terms, independent to the location of the particles, the complexity of the steps 2, 4 and 5 can be reduced.

\begin{table}
\begin{center}
   \begin{tabular}{| c | c |}
     \hline
     Step number & Complexity \\ \hline
     1 & $dN$ \\ \hline
     2 & $d^3 (2^D)^{\kappa}$ \\ \hline
     3 & $d^2 (2^D)^{\kappa}$ \\ \hline
     4 & $d^2 n_I (2^D)^\kappa$ \\ \hline
     5 & $d^3 (2^D)^\kappa$ \\     \hline
     6 & $d^2 (2^D)^{\kappa}$ \\ \hline
     7 & $dN$ \\     \hline
     $f_{\mathcal{N}}$ & $N^2 (2^D)^{-\kappa}$\\
     \hline
   \end{tabular}
 \end{center}
\caption{Complexity of each step of Algorithm~\ref{algo2} (for the most expensive level) and of the near interactions term ; $n_I$ denotes the maximum number of elements in the
interaction list of a box (4 in 1D, 40 in 2D and 316 in 3D).}
\label{tab22}
\end{table}

\section{Optimizing the overall complexity}
\label{sec:optimization}

\subsection{Precomputations for the M2M and L2L steps}

Recall that when using an EIM approximation, there holds $\Delta = B^{-t}\Gamma^{-1}$, where the matrices $B$ and $\Gamma$ are triangular and constructed in a precomputation step.
In practice, the condition number of the matrix $\Gamma$ worsens as the EIM approximation gets more accurate. It is well known in linear algebra that when a matrix $\Gamma$
is ill-conditioned, a direct resolution of the linear system $\Gamma x = b$ for some $b\in\mathbb{R}^d$ provides a result much less polluted by numerical errors than first inverting the matrix $\Gamma$ and then computing
the matrix vector product $\Gamma^{-1} b$. Since the matrices $\Gamma$ and $B$ are triangular, computing $\Delta b = y$ by first solving $\Gamma x = b$ and then solving $B^t y = x$
is of the same complexity $d^2$ as the direct matrix-vector product.

Consider step~2 in Algorithm~\ref{algo2} at level $k$:
\begin{equation*}
{W}_m^{I^{k}}=\sum_{p,q=1}^{d^{k+1}}\Delta_{p,q}^{I_0^{k+1}, J_0^{k+1}}\sum_{J^{k+1}\in\mathcal{C}(I^{k})}K(x_m^{I_0^{k}}+c_{I^{k}}-c_{J^{k+1}},y_p^{J_0^{k+1}})W_q^{J^{k+1}}.
\end{equation*}
We see that ${W}_m^{I^{k}}=\sum_{J^{k+1}\in\mathcal{C}(I^{k})} K^{I^k, J^{k+1}}\Delta^{I_0^{k+1}, J_0^{k+1}} W^{J^{k+1}}$, where $K^{I^k, J^{k+1}}_{m,p} = K(x_m^{I_0^{k}}+c_{I^{k}}-c_{J^{k+1}},y_p^{J_0^{k+1}})$,
$1\leq m\leq d^k$, $1\leq p\leq d^{k+1}$. The matrix $K^{I^k, J^{k+1}}_{M2M}:=K^{I^k, J^{k+1}}\Delta^{I_0^{k+1}, J_0^{k+1}}$ can be precomputed as $\left((\Delta^{I_0^{k+1}, J_0^{k+1}})^t(K^{I^k, J^{k+1}})^t\right)^t$, where
each matrix-vector product for each column of $(K^{I^k, J^{k+1}})^t$ is computed solving linear systems involving $B$ and $\Gamma$ matrices to preserve numerical accuracy, as explained in the
previous paragraph.
Notice that $K^{I^k, J^{k+1}}_{M2M}$ only depends on the relative position of $I^k$ and $J^{k+1}$. Since $J^{k+1}$ is a child of $I^k$, there are actually only $2^D$ different
operators $K^{I^k, J^{k+1}}_{M2M}$, where $D$ denotes the dimension of the space containing the particles. We denote them $K^{i}_{M2M}$, $1\leq i\leq 2^D$, and step 2 in
Algorithm~\ref{algo2} at level $k$ is reduced to
\begin{equation}
\label{eq:newM2M}
{W}_m^{I^{k}}=\sum_{i=1}^{2^D}K^{i}_{M2M}W_q^{J^{k+1}_i},
\end{equation}
where $J^{k+1}_i$ denotes the $i^{\rm th}$ child of $I^{k}$. This step is now of complexity $d^2 2^D (2^D)^k$ at level $k$, and the leading order of step 2 in Algorithm~\ref{algo2}
is $d^2 2^D (2^D)^{\kappa-1}=d^2 (2^D)^{\kappa}$.

The same optimization is done for step~5 in Algorithm~\ref{algo2}, reducing the complexity of this step to the order $d^2 (2^D)^{\kappa}$.

\subsection{Optimization of the M2L step}
\label{sec:compM2L}

It is well known that the M2L step is the most expensive one. Numerous efforts have been made to reduce the execution time of this step without
degrading the accuracy of the overall result, see in particular~\cite{DBLP:M2L}.
We use the optimization denoted SArcmp in~\cite{DBLP:M2L}.
It consists of a bivariate ACA applied to the collection of all the M2L operators at each level, followed by a truncated SVD.
The bivariate ACA and truncated SVD is then applied for each M2L operator at each level to increase further the compression.
The procedure is detailed for symmetric kernels for completeness of the presentation.

Consider a level $k$ and the collection of all the M2L operators at this level.
Consider a level $k$, a box $I^k$ of the tree at this level $k$, and a box $J^k$ in the interaction list of $I^k$. The vector $\delta:=c_{J^k}-c_{I^k}$ only depends on the relative position
of the boxes $I^k$ and $J^k$. Therefore, there are only $n_I = 7^D-3^D$ different M2L matrices $K^\delta_{i,j}:=K(x_{i}^{J_0^k},y_j^{J_0^k}+\delta)$ per level (we recall that $D$
denotes the dimension of the space containing the particles), see~\eqref{eq:M2Lstep}.
We suppose that EIM algorithms have been carried-out for a certain error bound $\epsilon$, leading to the choice of $d_k$ interpolation points.
The M2L matrices are organized in one large matrix as follows $K_{\rm{fat}}:=[K^{\delta_1} K^{\delta_2} ... K^{\delta_{n_I}}]$, where $K_{\rm{fat}}\in\mathbb{R}^{d_k\times n_I d_k}$.
Two rectangular matrices $U,V$ are computed by a bivariate ACA such that $|UV|_F\leq \epsilon|K_{\rm{fat}}|_F$,
where $|\cdot|_F$ denotes the Frobenius norm. Then, two QR decompositions are computed, such that $U=Q_U R_U$ and $(V)^t=Q_V R_V$,
and a SVD is computed on $R_U(R_V)^t$. Denote $s_i$, $1\leq i\leq d_k$, the singular values of $R_U(R_V)^t$ in decreasing order, and $\hat{U}$, $\hat{V}$
the unitary matrices such that $R_U(R_V)^t=\hat{U}{\rm diag}(s)\hat{V}$.
The two following filters are applied: (i) keep the first $j$ singular values such that $\frac{s_j}{s_0}\leq\epsilon\leq \frac{s_{j-1}}{s_0}$, (ii) keep the first $j$ singular values such that 
$\frac{\sum_{i=1}^{j}s_i}{\sum_{i=1}^{d^k}s_i}\geq1-\epsilon\geq \frac{\sum_{i=1}^{j-1}s_i}{\sum_{i=1}^{d^k}s_i}$.
We denote $r_k$ the number of kept singular values. Denote $\hat{U}_{r_k}$ and $\hat{V}_{r_k}$ the first $r_k$ columns of respectively $\hat{U}$ and $\hat{V}$.
Following~\cite{bbFMM}, we can show that the M2L step can be reduced to size $r_k\times r_k$ (instead of size $d_k\times d_k$) matrix-vector products
involving the matrices $C^{\delta}:=\hat{U}_{r_k}^T K^{\delta}\hat{U}_{r_k}$, plus low complexity pre- and post-processing steps.

Then, the same compression is applied locally for each compressed M2L operator $C^{\delta}$. This time, the bivariate ACA is applied on each $C^{\delta}$, producing the rectangular matrices $U$ and $V$.
The QR decompositions, SVD and filters are applied in the same fashion, selecting the singular values $s_k$.
The rectangular matrices $\hat{U}_{s_k}$ and $\hat{V}_{s_k}$ are constructed, and
a diagonal matrix $S_{s_k}$, whose diagonal entries are the square roots of the selected singular values, is constructed.
Then, we construct $\tilde{U}_{s_k}=Q_U\hat{U}_{s_k}S_{s_k}\in\mathbb{R}^{r_k\times s_k}$ and $\tilde{V}_{s_k}=S\hat{V}_{s_k}Q_V^t\in\mathbb{R}^{s_k\times r_k}$,
so that $C^{\delta}$ is approximated by $\tilde{U}_{s_k} \tilde{V}_{s_k}$. Thus, size $r_k\times r_k$ matrix-vector products involving the matrices $C^{\delta}$
are replaced by two $r_k\times s_k$ matrix-vector products involving the matrices $\tilde{U}_{s_k}$ and $\tilde{V}_{s_k}$.

\section{Numerical experiments in 3D}
\label{sec:num}

\subsection{Some elements on the implementation}

The EIFMM has been implemented in C++, using the open source library scalfmm, see~\cite{scalfmm}.
From~\cite{benchfmm}, scalfmm is one of the fastest opensource library for FMM computations.
We will mainly compare our implementation of the EIFMM to the Chebyshev interpolation-based implementation in scalfmm, which follows the Black-Box FMM~\cite{bbFMM}.
In our code, the EIM precomputation and the SArcmp precomputation steps are performed in C++ using BLAS and LAPACK, and stored in binary files.
These precomputations are then read by the modified scalfmm library, using also BLAS for the basic linear algebra operations.

\subsection{Test-cases}

\begin{figure}[h!]
   \begin{minipage}[c]{.32\linewidth}
\includegraphics[width=\textwidth]{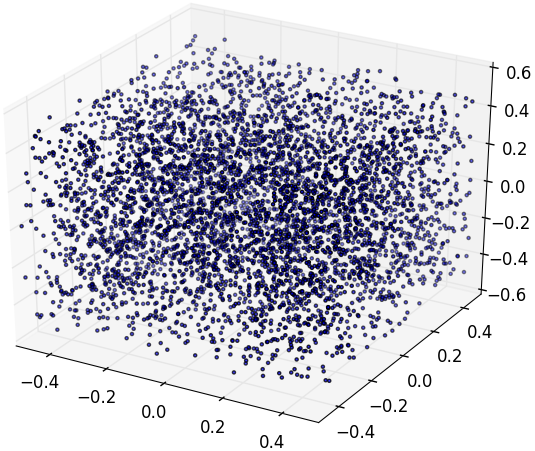}
   \end{minipage} \hfill
   \begin{minipage}[c]{.32\linewidth}
\includegraphics[width=\textwidth]{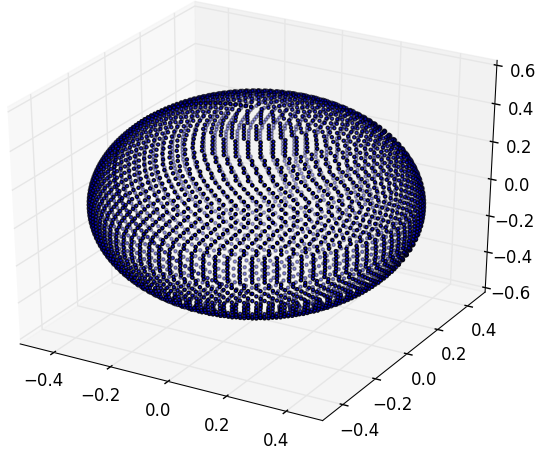}
   \end{minipage}
   \begin{minipage}[c]{.32\linewidth}
\includegraphics[width=\textwidth]{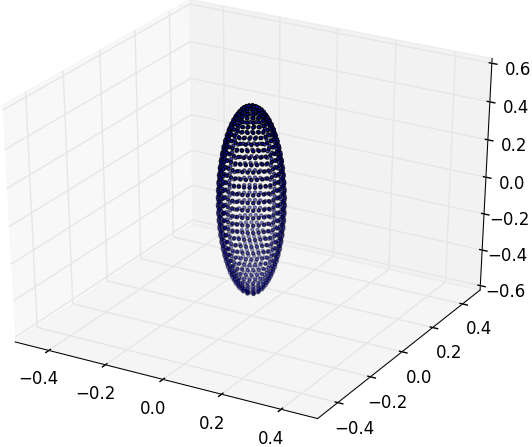}
   \end{minipage}
 \caption{Sets of points considered in the numerical applications, from left to right: random ($5000$ points), sphere ($10100$ points), and ellipsoid ($1944$ points).
The values for the number of points have been chosen for representation reasons, higher values will be considered in the numerical experiments.}
\label{fig:pointsets}
\end{figure}

We consider sets of points included in $\Omega=(-0.5, 0.5)^3$, randomly taken following (i) a uniform law in $\Omega$, (ii) a uniform law on the sphere of radius $0.5$ and
(iii) a modification of the set (ii) in an ellipsoid, see Figure~\ref{fig:pointsets}.

When comparing runtimes, we have to make sure that the Black-Box FMM is used in a optimal way. There are two parameters to choose: the number of interpolation points
per dimension and the accuracy $\epsilon_c$ of the compression. First, we chose $3$ interpolation points per dimension
and compute the relative error with $\epsilon_c=10^{-15}$. Then, if this relative error is smaller than the desired one, we search by dichotomy
the largest value of the $\epsilon_c$ so that the relative error is smaller or equal to the desired one. Otherwise, we increase the number of interpolation points
per dimension by one and repeat the process. That way, the runtimes comparison is fair in the sense that Block-Box FMM is not artificially slowed down by
a too large number of interpolation points together with a too large value for $\epsilon_c$. Notice that this trick is necessary since we do not know in advance
how the Chebyshev interpolation will perform in practice for a new kernel, while for the EIFMM, the error estimate automatically selects the right number of
interpolation points.

\subsection{Performance for the Laplace kernel}

In this section, we compare runtimes of the FMM summation using the Black-Box FMM and the EIFMM.
The computations are carried out in sequential, on a laptop with 4~GB of RAM and a Intel(R) Core(TM) i5 CPU M 560 @ 2.67GHz.
For the Black-Box FMM, we use
the so-called symmetric M2L compression, based on local SVD for each M2L operator. The symmetric M2L compression
exploit the symmetry of the Chebyshev nodes to express the M2L step as matrix-matrix products, and even if the complexity of this step is not improved, 
runtimes are drastically improved due to better reuse of memory cache by the processor. Proposed in~\cite{DBLP:M2L}, the symmetric compression is currently, to the author's knowledge,
the most efficient compression strategy for the Black-Box FMM.
Since the approximation~\eqref{eq:approxFMM} used in EIFMM has four summations instead of two for Black-Box FMM, the steps 3 and 6 in Algorithm~\ref{algo2} for EIFMM
are not present in Black-Box FMM.

\begin{figure}[h!]
\begin{center}
\includegraphics[width=0.5\textwidth]{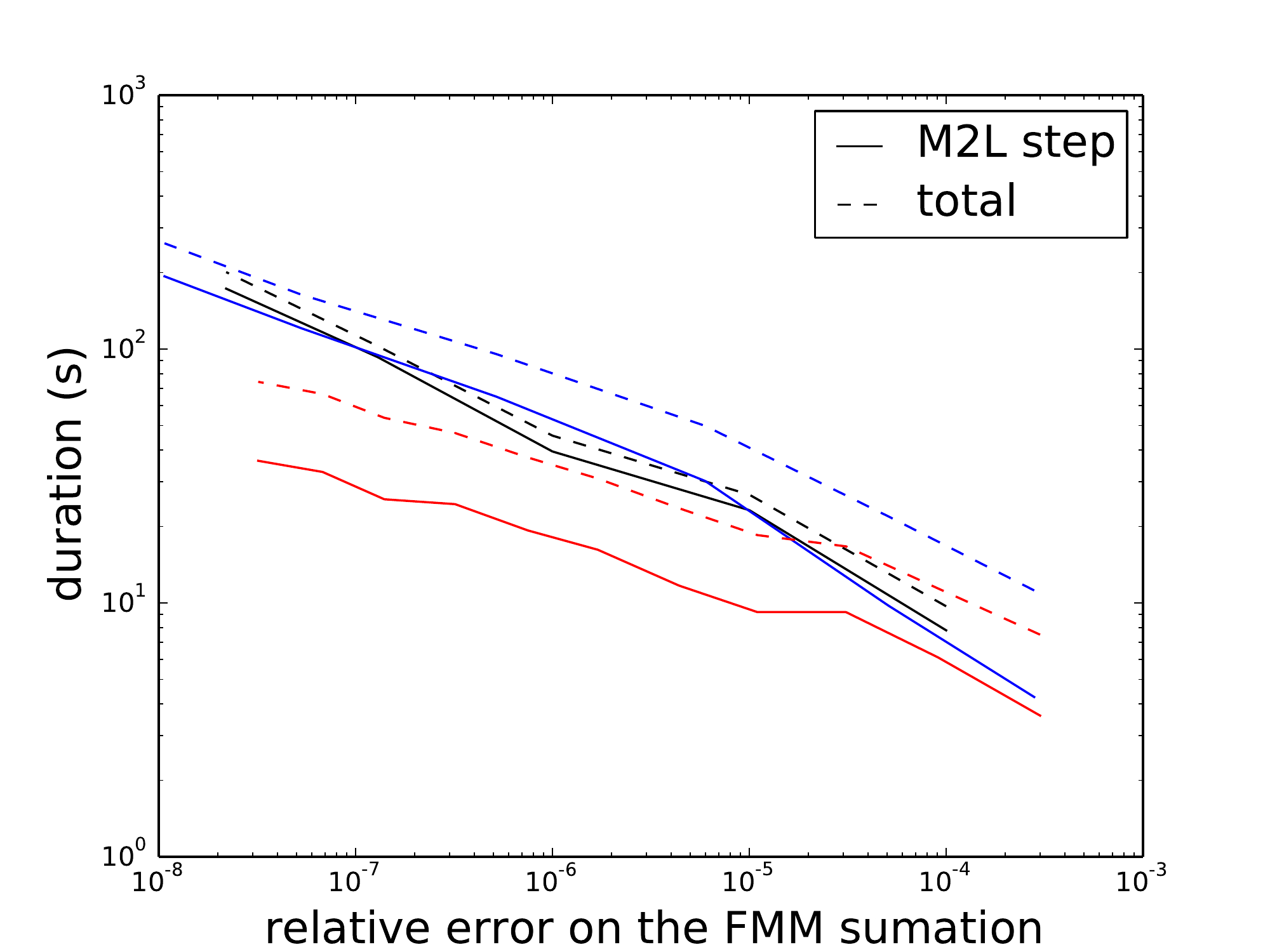}
\end{center}
 \caption{Execution time with respect to the relative error on the FMM summation for the Laplace kernel;
black: symmetric Black-Box FMM, blue: EIFMM, red: fastest spherical algorithm available in scalfmm.}
\label{fig:1sr}
\end{figure}

In Figure~\ref{fig:1sr} are represented execution time with respect to the relative error on the FMM summation for the Laplace kernel,
for the symmetric Black-Box FMM, the EIFMM and the fastest algorithm available in scalfmm.
This latter is dedicated to the Laplace kernel, explaining why it performs better than the other two.
The test-case corresponds to $10^6$ random points taken in the unit cube, organized in a 6-level tree, see Figure~\ref{fig:pointsets}, left picture.
The execution times are better for Black-Box FMM at low accuracy, and equivalent between EIFMM and Black-Box FMM at high accuracy.
Notice that since the Laplace kernel is homogeneous, the difficulty of the approximation problem does not change from one level to another in the tree.
As a consequence, the variable-order capability of our algorithm is not exploited for this kernel. In what follows, inhomogeneous kernels are considered and gains
on runtimes are observed. Precomputation times will be investigated in the following section as well.

\subsection{Performance for some inhomogeneous kernels}
\label{sec:num_other}

The advantages of the EIFMM mainly apply to inhomogeneous kernels, for which the error estimation allows the selection of an optimal number of interpolation
points at each level of the tree by the EIM.

\subsubsection{1-million-point test-cases}

The three test-cases presented in Figure~\ref{fig:pointsets} are considered: $10^6$ random points in the unit cube with a 6-level tree (denoted ''cube``),
$992,\!046$ points at the surface of a sphere and an ellipsoid, with respectively 8-level and 9-level trees (denoted respectively ''sphere`` and ''ellipsoid``).

\begin{figure}[h!]
   \begin{minipage}[c]{.5\linewidth}
\includegraphics[width=\textwidth]{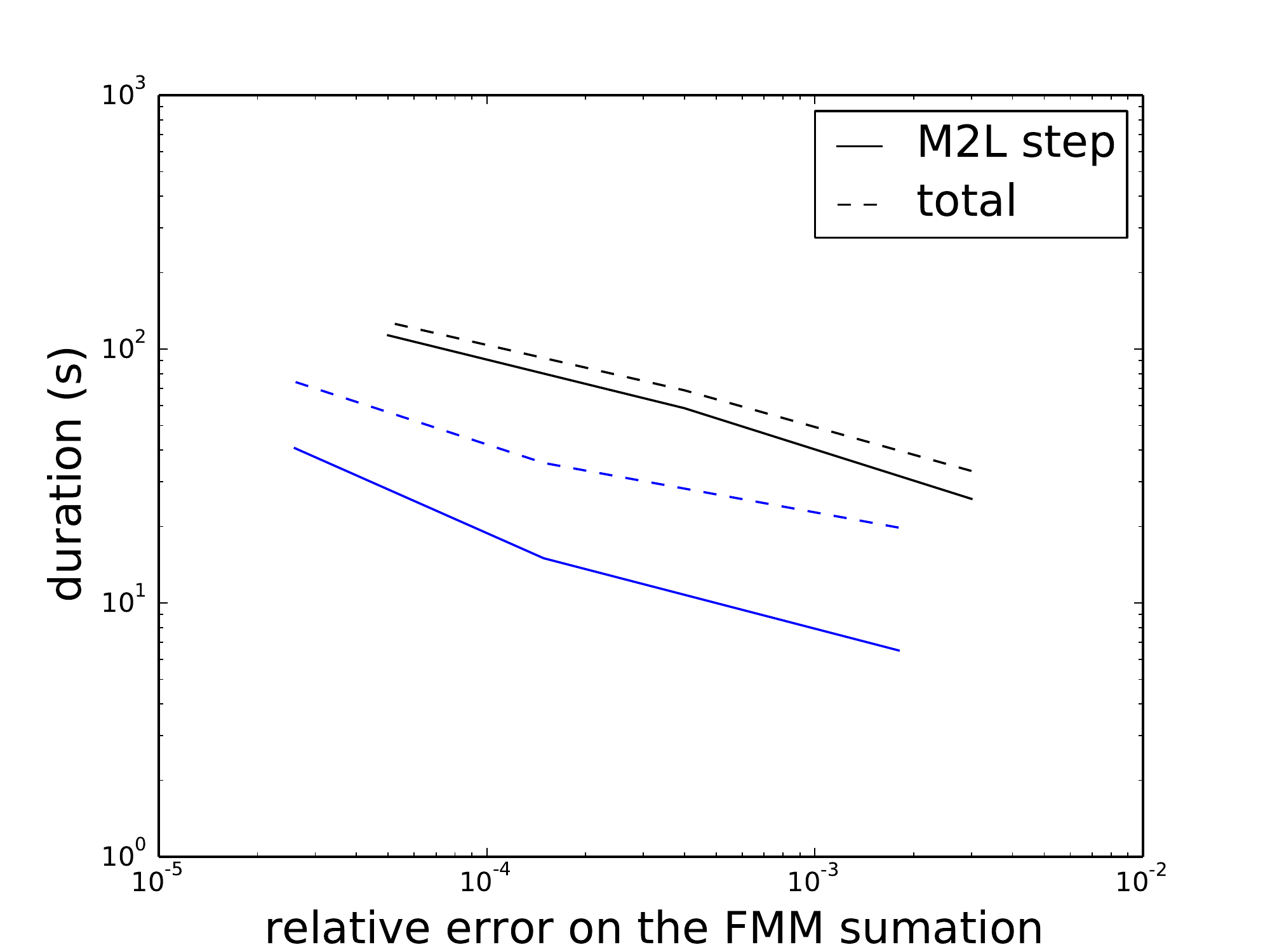}
   \end{minipage} \hfill
   \begin{minipage}[c]{.5\linewidth}
\includegraphics[width=\textwidth]{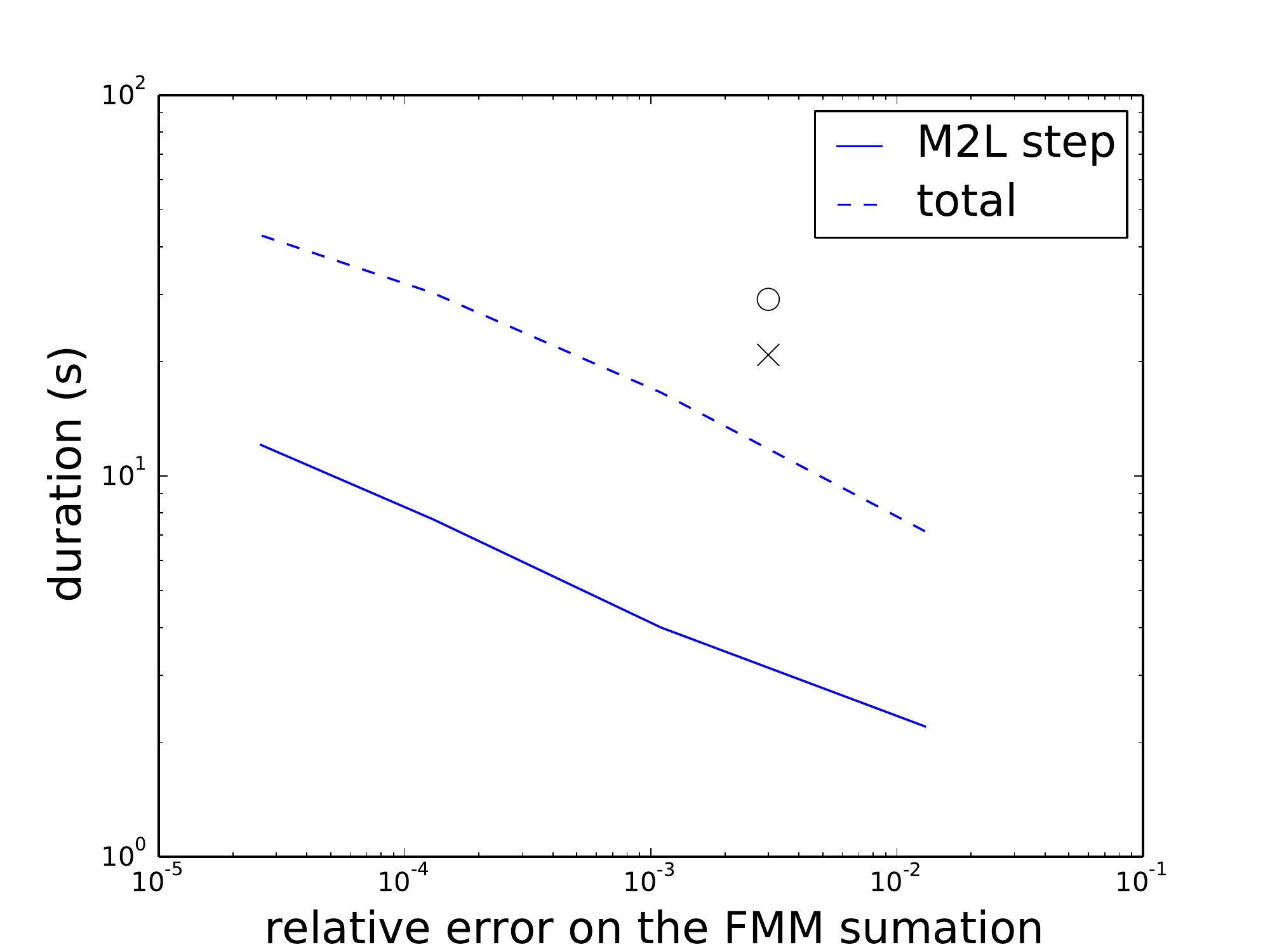}
   \end{minipage}
\begin{center}
\includegraphics[width=0.5\textwidth]{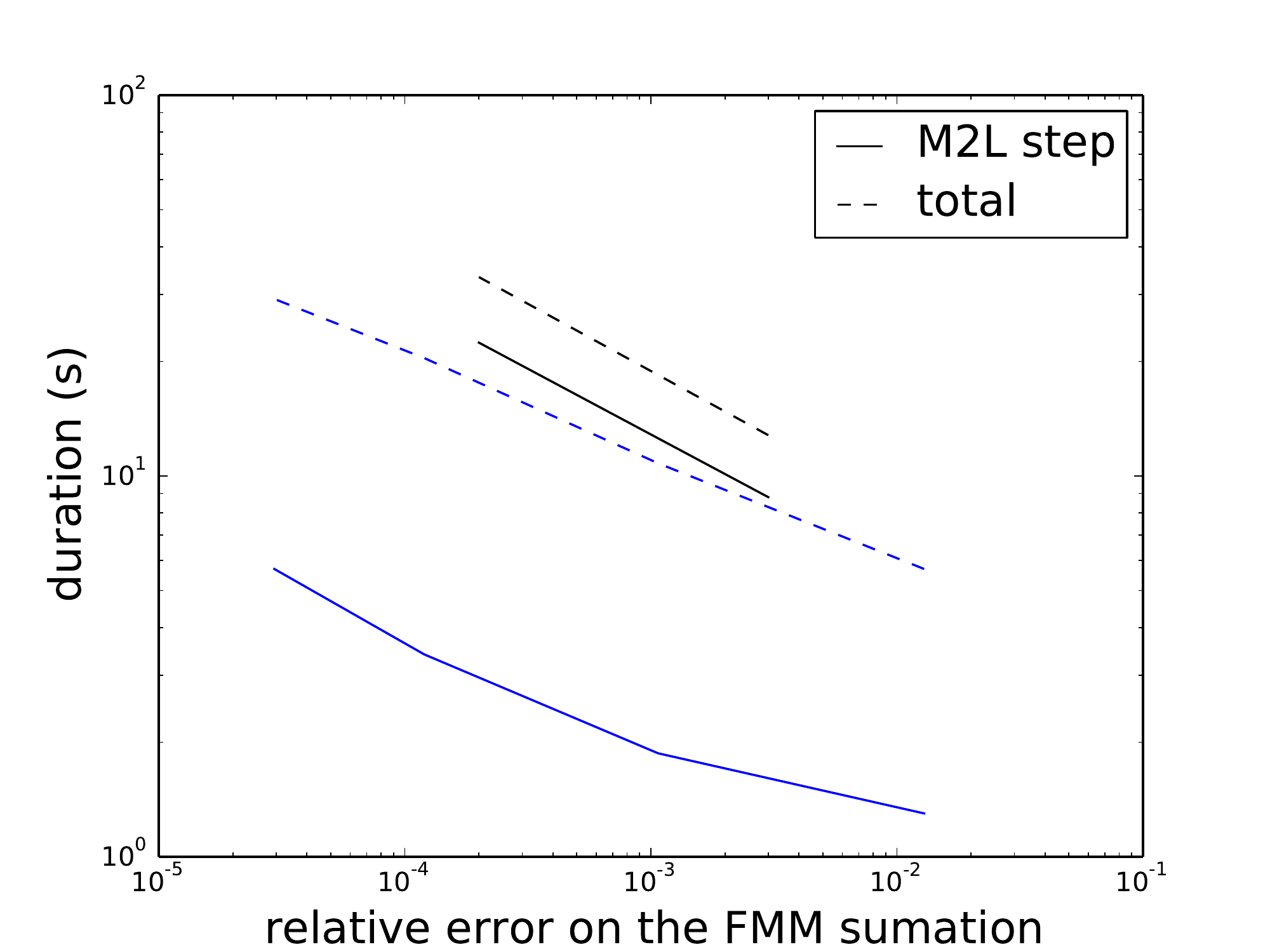}
\end{center}
 \caption{Execution time of the long-range FMM accelerated terms with respect to the relative error on the FMM summation for the kernel $\frac{\cos{20r}}{r}$, black: symmetric Black-Box FMM, blue: EIFMM.
From top to bottom and left to right: cube, sphere and ellipsoid test-cases.}
\label{fig:cos}
\end{figure}

\begin{figure}[h!]
   \begin{minipage}[c]{.5\linewidth}
\includegraphics[width=\textwidth]{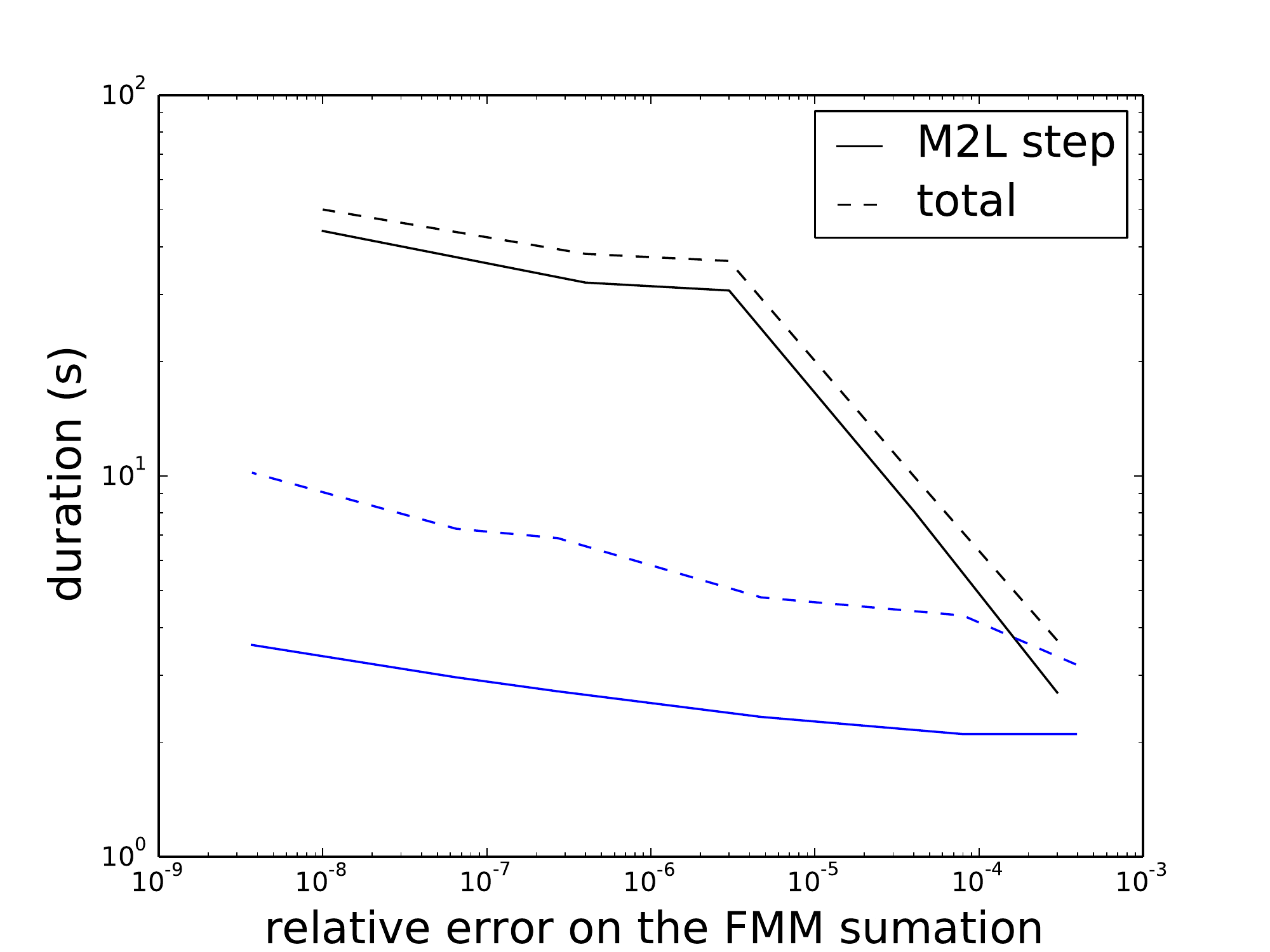}
   \end{minipage} \hfill
   \begin{minipage}[c]{.5\linewidth}
\includegraphics[width=\textwidth]{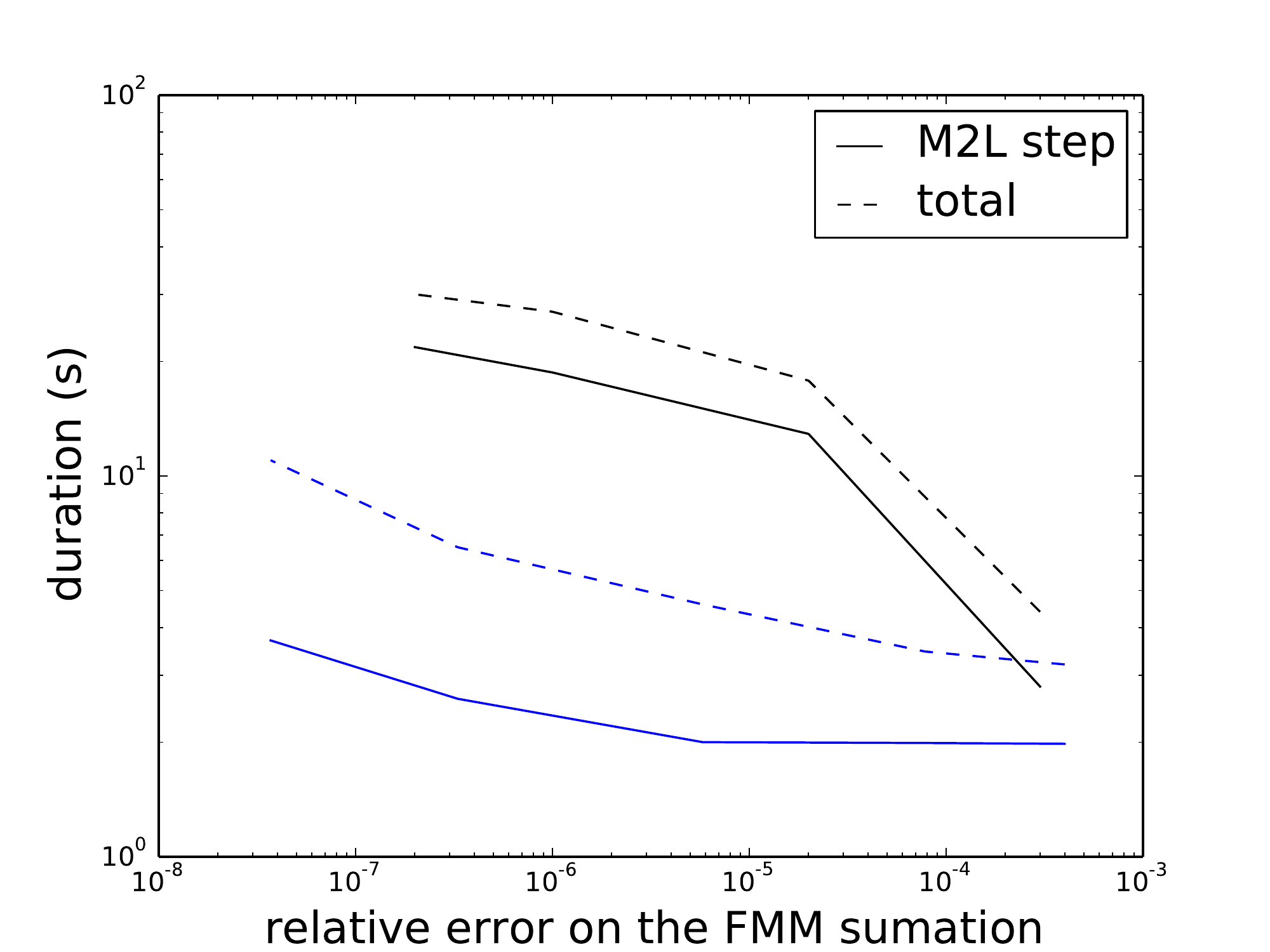}
   \end{minipage}
\begin{center}
\includegraphics[width=0.5\textwidth]{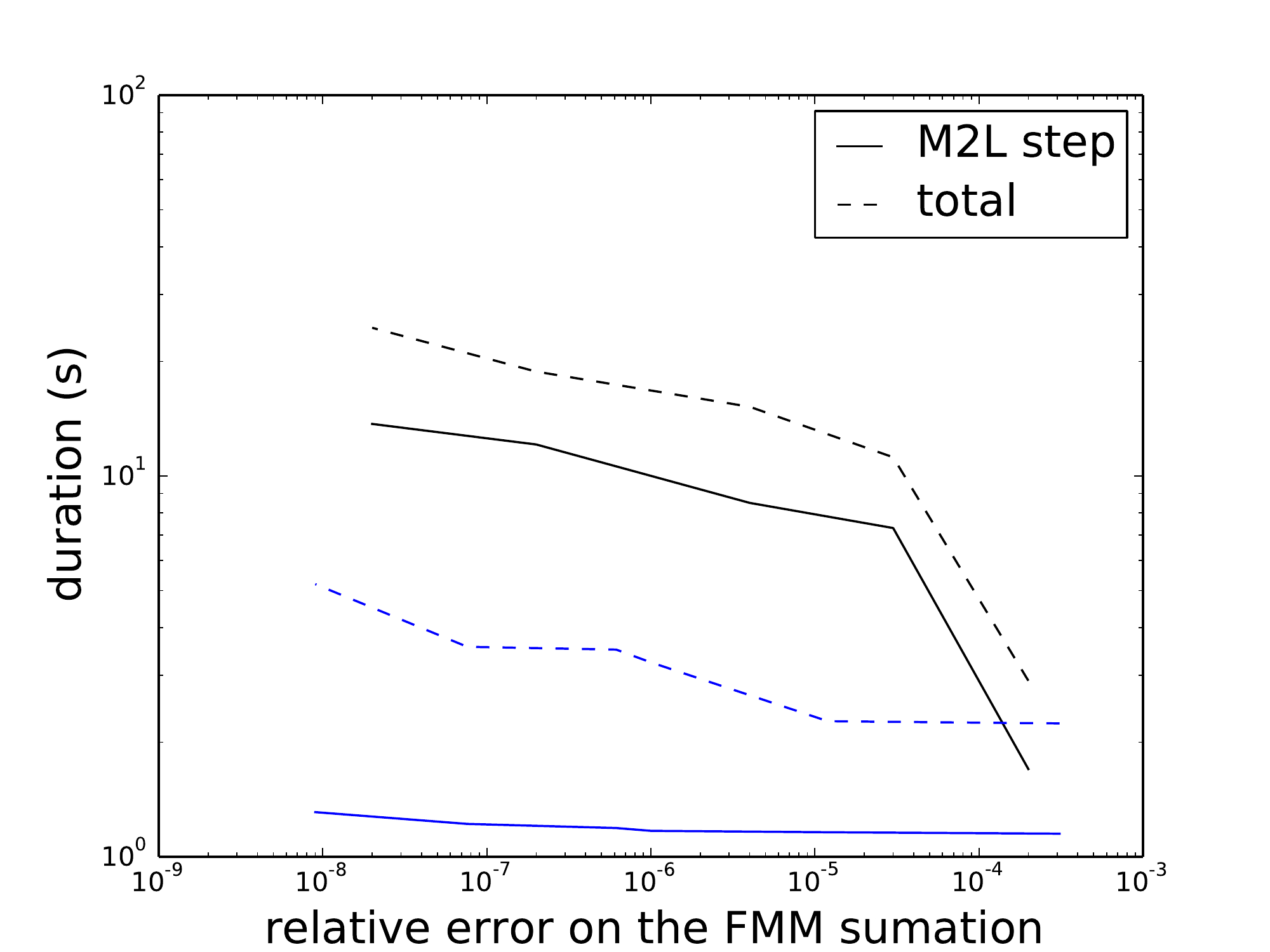}
\end{center}
 \caption{Execution time of the long-range FMM accelerated terms with respect to the relative error on the FMM summation for the kernel $e^{-r^2}$, black: symmetric Black-Box FMM, blue: EIFMM.
From top to bottom and left to right: cube, sphere and ellipsoid test-cases.}
\label{fig:e-r2}
\end{figure}

\begin{figure}[h!]
   \begin{minipage}[c]{.5\linewidth}
\includegraphics[width=\textwidth]{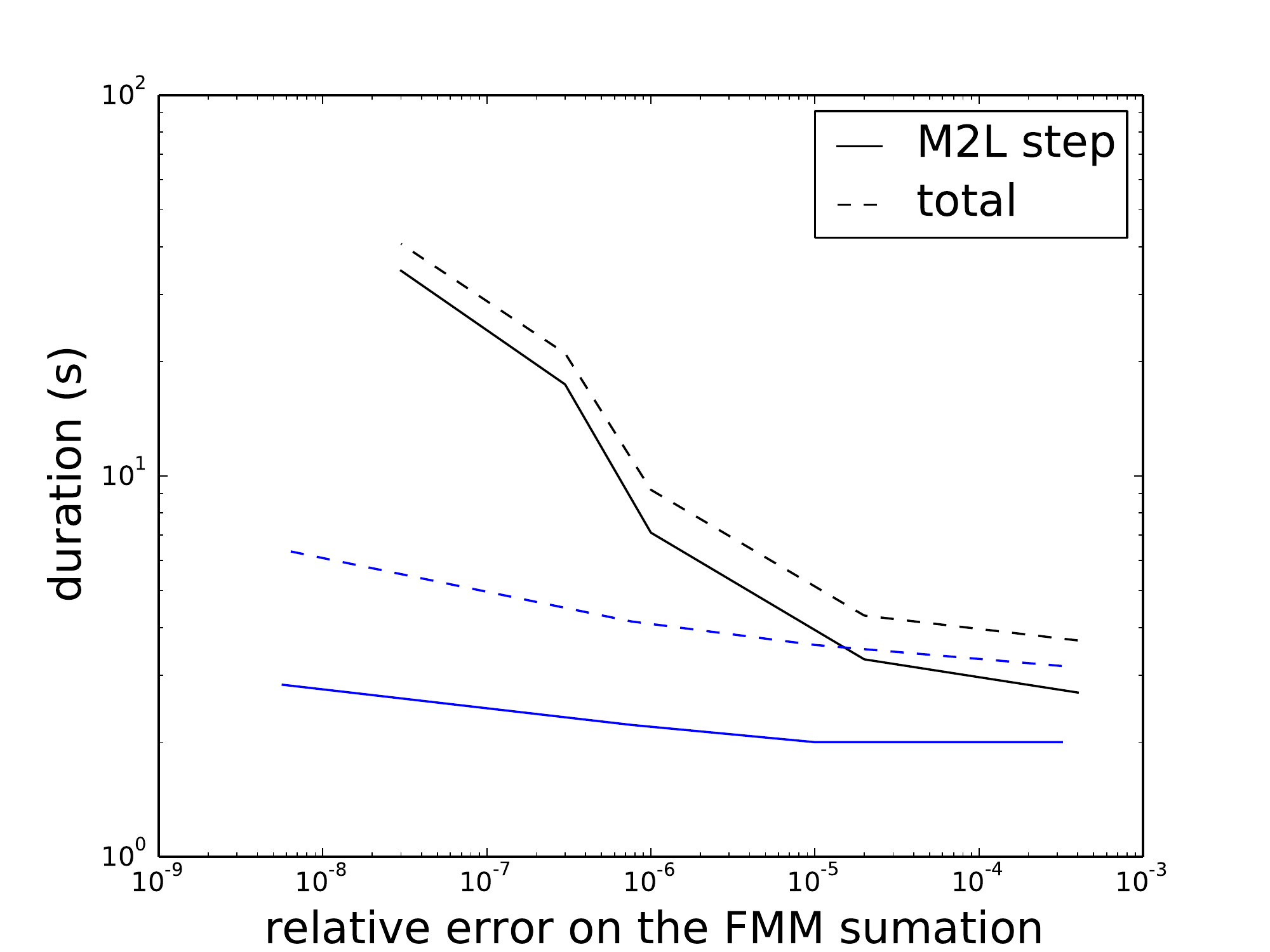}
   \end{minipage} \hfill
   \begin{minipage}[c]{.5\linewidth}
\includegraphics[width=\textwidth]{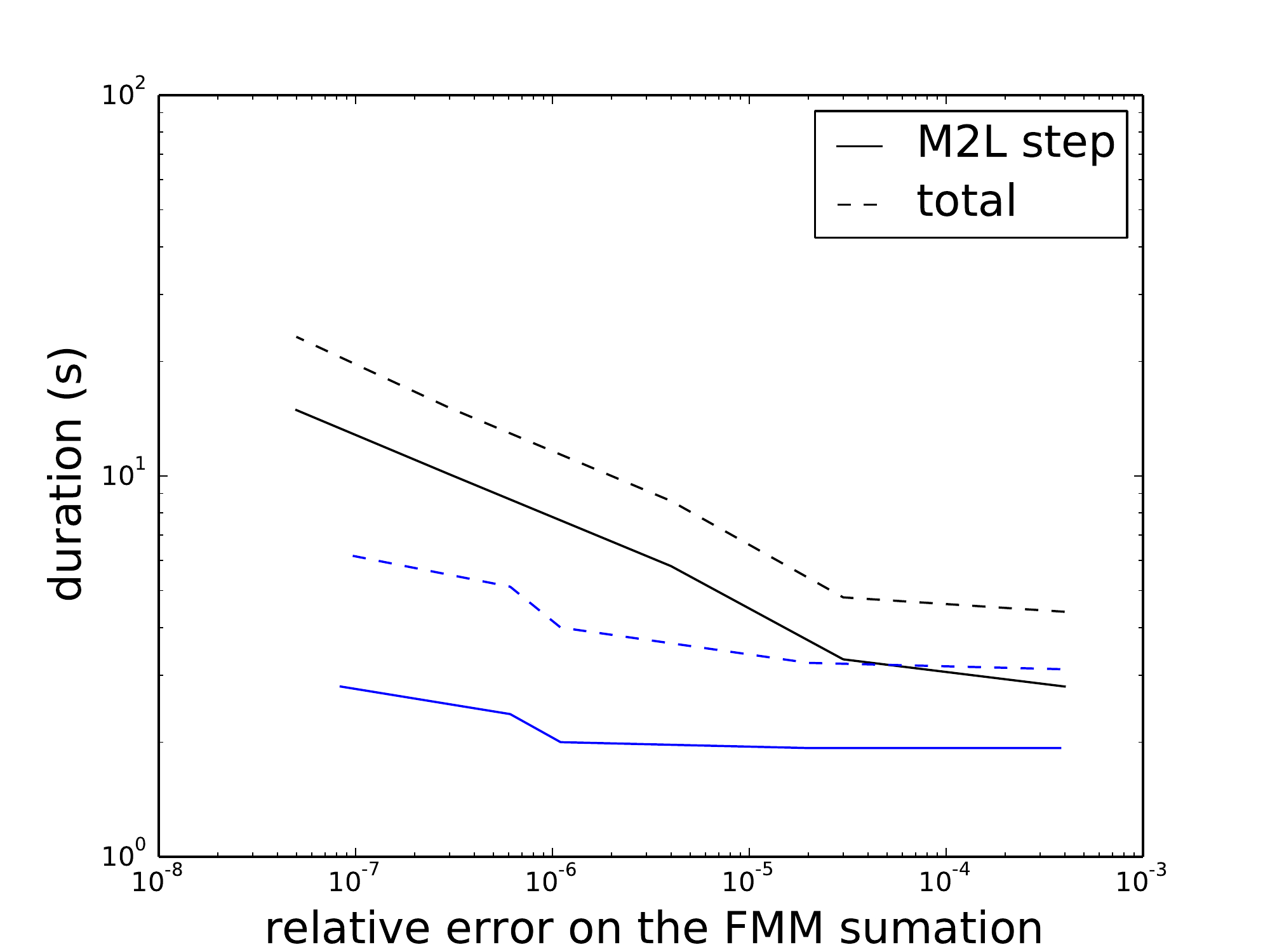}
   \end{minipage}
\begin{center}
\includegraphics[width=0.5\textwidth]{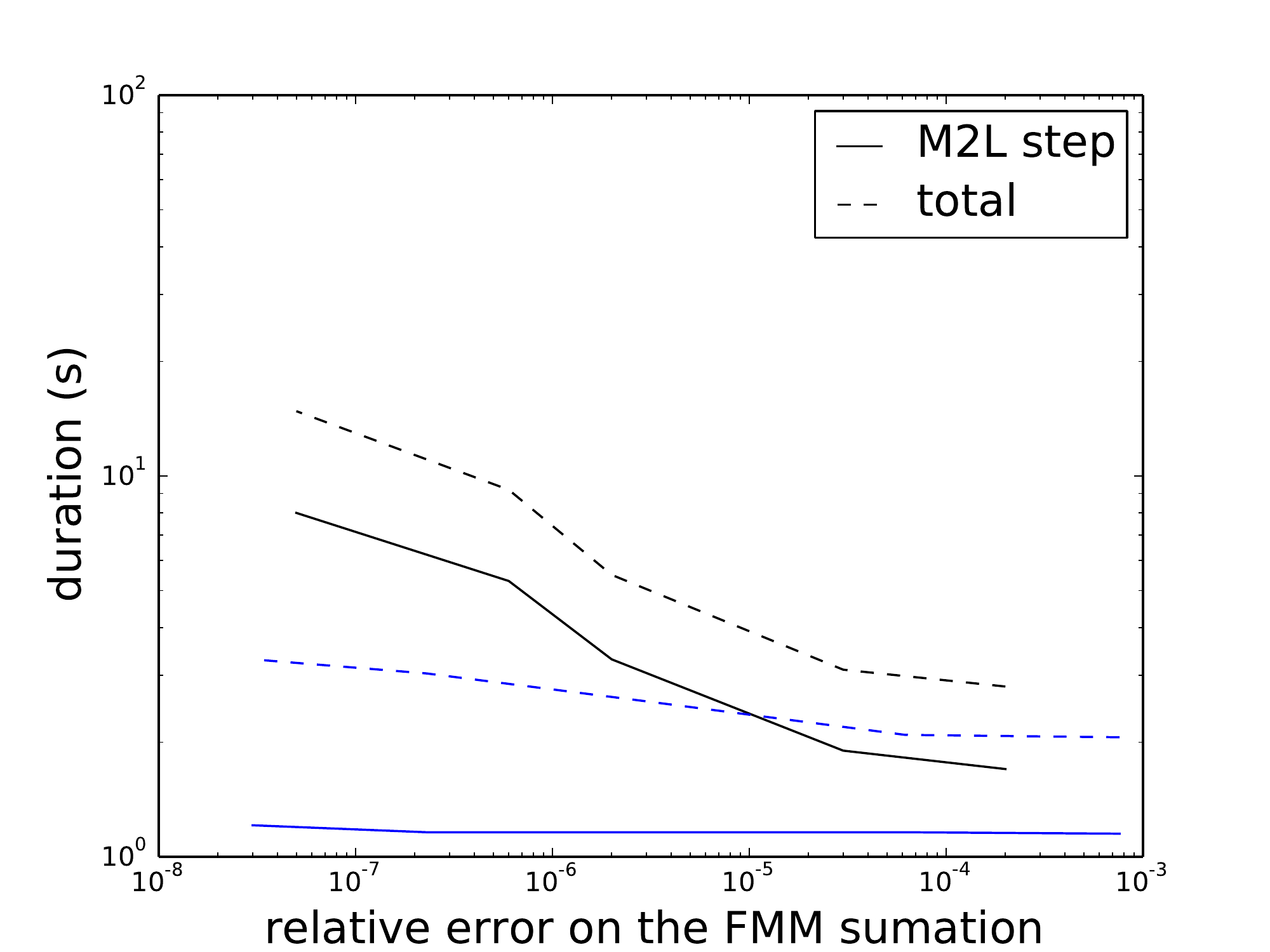}
\end{center}
 \caption{Execution time of the long-range FMM accelerated terms with respect to the relative error on the FMM summation for the kernel $\sqrt{r^2+1}$, black: symmetric Black-Box FMM, blue: EIFMM.
From top to bottom and left to right: cube, sphere and ellipsoid test-cases.}
\label{fig:sqrtr2+1}
\end{figure}

Figures~\ref{fig:cos},~\ref{fig:e-r2} and~\ref{fig:sqrtr2+1} present execution times of the long-range FMM accelerated terms with respect to the relative error on the FMM summation,
respectively for the kernels $\frac{\cos{20r}}{r}$, $e^{-r^2}$ and $\sqrt{r^2+1}$, for the three test-cases (cube, sphere and ellipsoid).
The reported runtimes do not contain the times required to build the tree, nor the EIM, M2L compression, M2M and L2L precomputations.
For the EIFMM, we use the SArcmp compression for the M2L operators, as described in Section~\ref{sec:compM2L}.
Due to smaller interpolation formulae, EIFMM enables memory savings with respect to Black-Box FMM. 
In the figures, missing data correspond to cases requiring more memory than available with our computer, for which EIFMM can be used, but not Blakc-Box FMM.
For the three kernels, and especially the last two, the EIFMM is much faster than the Black-Box FMM, because
the error estimation and the greedy procedure in EIFMM allows the selection of an optimal number of interpolation points at each level of the tree.
For instance, for the kernel $e^{-r^2}$ on the ellipsoid and a relative error of $10^{-6}$, the Black-Box FMM needs 6 interpolation points per dimension (i.e. 216
interpolation points), whereas the EIFMM needs respectively 167, 74, 36, 21, 13, 10, 9, and 4 points at levels 2 to 9 of the tree. As we go deeper in the tree, the
interpolation problems become easier, and the EIM selects fewer points to produce an interpolation of the same quality.

Table~\ref{tab_precomp_times} provides precomputations and executions times for the kernels $\frac{\cos(20r)}{r}$, $e^{-r^2}$ and $\sqrt{r^2+1}$ on the cube test-case.
Notice that in order not to artificially disadvantage Black-Box FMM, we report values with EIFMM set to a better relative accuracy and we neglect the precomputation times
of Black-Box FMM (which are short anyway). In the reported cases, EIFMM becomes advantageous (i.e. the EIFMM precomputations times are compensated by faster 
long-range FMM terms computation time) if one has to compute more than 2 to 18 FMM summations.

\begin{table}[h]
\begin{center}
\begin{tabular}{|c|c|c|c|c|c|c|}
\hline
kernel                             & \multicolumn{2}{c|}{$\frac{\cos(20r)}{r}$} & \multicolumn{2}{c|}{$e^{-r^2}$}       & \multicolumn{2}{c|}{$\sqrt{r^2+1}$}   \\ \hline
algorithm                          & EIFMM                & BBFMM               & EIFMM             & BBFMM             & EIFMM             & BBFMM             \\ \hline
relative accuracy                  & $2\times 10^{-4}$    & $4\times 10^{-4}$   & $2\times 10^{-6}$ & $3\times 10^{-6}$ & $1\times 10^{-5}$ & $2\times 10^{-5}$ \\ \hline
EIM precomp time (s)               & 73                   & --                  & 30                & --                & 10                & --                \\ \hline
compression time (s)               & 177                  & --                  & 33                & --                & 1.1               & --                \\ \hline
tree construction time (s)         & \multicolumn{2}{c|}{3}                     & \multicolumn{2}{c|}{3}                & \multicolumn{2}{c|}{3}                \\ \hline
long-range FMM terms comp time (s) & 52                   & 69                  & 5.4               & 37                & 3.5               & 4.3               \\ \hline
short-range terms comp time(s)     & \multicolumn{2}{c|}{32}                    & \multicolumn{2}{c|}{30}               & \multicolumn{2}{c|}{16}               \\ \hline
minimum EIFMM calls                & \multicolumn{2}{c|}{18}                    & \multicolumn{2}{c|}{2}                & \multicolumn{2}{c|}{14}               \\ \hline
\end{tabular}
 \end{center}
\caption{Precomputations and executions times for the kernels $\frac{\cos(20r)}{r}$, $e^{-r^2}$ and $\sqrt{r^2+1}$ on the cube test-case.
The last line indicates the minimum FMM summations to compute for EIFMM to be faster than Black-Box FMM (BBFMM), taking all precomputations times into account.}
\label{tab_precomp_times}
\end{table}

\subsubsection{100-million-point test-case}

We consider $10^8$ random points in the unit cube (Figure~\ref{fig:pointsets}, left) in a 9-level tree, with the $e^{-r^2}$ kernel.
The computations are carried-out in sequential on a Intel(R) Xeon(R) CPU E5-4620 v2 @ 2.60GHz and 128~GB of RAM.

Table~\ref{100mill} presents the relative accuracy and the execution times of all the steps for two simulations using EIFMM and Black-Box FMM.
For a better accuracy, the EIFMM is faster than the Black-Box FMM, including the precomputation times (increasing the accuracy with Black-Box FMM would require more than
128~GB of RAM). In this case, EIFMM is interesting even for a single FMM computation.
Notice that the precomputations for EIFMM are negligible for such large cases: the precomputation complexity is linear with respect to the number of levels in the tree
(therefore logarithmic with respect to $N$).
Moreover, for inhomogeneous kernels, the approximation problem for an additional level is easier than the previous approximation problems, leading to even faster additional EIM 
and M2L compression precomputation times.

\begin{table}[h]
\begin{center}
\begin{tabular}{|c|c|c|}
\hline
algorithm                  & EIFMM             & Black-Box FMM     \\ \hline
relative accuracy          & $7\times 10^{-5}$ & $6\times 10^{-4}$ \\ \hline
EIM precomp time           & 22.8 s            & --                \\ \hline
compression time           & 1.3 s             & --                \\ \hline
tree construction time (s) & \multicolumn{2}{c|}{9 min 10 s}       \\ \hline
P2M (1)                       & 1 min 13 s        & 2 min 17 s        \\ \hline
M2M (2-3)                       & 31 s              & 49 s              \\ \hline
M2L (4)                       & 10 min 39 s       & 53 min 01 s       \\ \hline
L2L (5-6)                       & 41 s              & 59 s              \\ \hline
L2P (7)                       & 1 min 42 s        & 3 min 55 s        \\ \hline
total long-range FMM terms & 14 min 46 s       & 1 h 01 min 01 s   \\ \hline
short-range terms          & 19 min 36 s       & 19 min 36 s       \\ \hline
total                      & 43 min 56 s       &  1 h 29 min 47 s    \\ \hline
\end{tabular}
 \end{center}
\caption{Accuracy and execution times for the 100 million points test-case. P2M, M2M, M2L, L2L and L2P corresponds
respectively to the steps (1), (2-3), (4), (5-6) and (7) in Algorithm~\ref{algo2}.}
\label{100mill}
\end{table}

\section{Conclusion and outlook}
\label{sec:conclusion}

This work introduces a new Fast Multipole Method (FMM) using a low-rank approximation of the kernel based on the Empirical Interpolation Method (EIM), called
the Empirical Interpolation Fast Multipole Method (EIFMM).
The proposed multilevel algorithm is implemented in scalfmm, a FMM library written in C++.
The important feature of the EIFMM is a built-in error estimation of the interpolation error made by the low-rank approximation of the far-field behavior of the kernel.
As a consequence, the algorithm automatically selects the optimal number of interpolation points required to ensure a given accuracy for the result, 
leading to important execution time reduction for inhomogeneous kernels, for which the difficulty of the approximation varies from one level of the tree to another.

\section*{Acknowledgement}
The author wish to thank Centre National d'Etudes Spatiales for financial support through the TOSCA committee, 
Guillaume Sylvand (Airbus Group Innovations, INRIA) for drawing his attention to fast multipole methods, thoroughly reading the manuscript and for 
many fruitfull discussions, Tony Leli\`evre (CERMICS) and Alexandre Ern (CERMICS) for fruitful discussions, and Mathieu Ruffat (Diginext) for his C++ expertise.


\begin{thebibliography}{10}
\bibitem{wavelets}
B.~Alpert, G.~Beylkin, R.~Coifman, and V.~Rokhlin.
\newblock Wavelet-like bases for the fast solutions of second-kind integral
  equations.
\newblock {\em SIAM J. Sci. Comput.}, 14(1):159--184, 1993.

\bibitem{Barrault}
M.~Barrault, Y.~Maday, N.~C. Nguyen, and A.~T. Patera.
\newblock An 'empirical interpolation' method: application to efficient
  reduced-basis discretization of partial differential equations.
\newblock {\em Comptes Rendus Mathematique}, 339(9):667 -- 672, 2004.

\bibitem{ACA}
M.~Bebendorf.
\newblock Approximation of boundary element matrices.
\newblock {\em Numerische Mathematik}, 86(4):565--589, 2000.

\bibitem{Hmat}
M.~Bebendorf and S.~Rjasanow.
\newblock Adaptive low-rank approximation of collocation matrices.
\newblock {\em Computing}, 70(1):1--24, 2003.

\bibitem{bebendorf}
M.~Bebendorf and R.~Venn.
\newblock Constructing nested bases approximations from the entries of
  non-local operators.
\newblock {\em Numerische Mathematik}, 121(4):609--635, 2012.

\bibitem{borm}
S.~B\"orm, M.~L\"ohndorf, and J.M. Melenk.
\newblock Approximation of integral operators by variable-order interpolation.
\newblock {\em Numerische Mathematik}, 99(4):605--643, 2005.

\bibitem{casenave_ACM}
F.~Casenave, A.~Ern, and T.~Leli\`evre.
\newblock A nonintrusive reduced basis method applied to aeroacoustic
  simulations.
\newblock {\em Advances in Computational Mathematics}, pages 1--26, 2014.

\bibitem{lowrank}
H.~Cheng, Z.~Gimbutas, P.G. Martinsson, and V.~Rokhlin.
\newblock On the compression of low-rank matrices.
\newblock {\em SIAM J. Sci. Comput.}, 26(4):1389--1404, 2005.

\bibitem{scalfmm}
O.~Coulaud, B.~Bramas, and C.~Piacibello.
\newblock Scalfmm, {C}++ {F}ast {M}ultipole {M}ethod {L}ibrary for {HPC}.
\newblock {\em http://scalfmm-public.gforge.inria.fr/doc/}.

\bibitem{cheby0}
A.~Dutt, M.~Gu, and V.~Rokhlin.
\newblock Fast algorithms for polynomial interpolation, integration, and
  differentiation.
\newblock {\em SIAM Journal on Numerical Analysis}, 33(5):1689--1711, 1996.

\bibitem{fft}
A.~Dutt and V.~Rokhlin.
\newblock Fast fourier transforms for nonequispaced data.
\newblock {\em SIAM Journal on Scientific Computing}, 14(6):1368--1393, 1993.

\bibitem{EIMeff3}
J.~L. {Eftang} and B.~{Stamm}.
\newblock Parameter multi-domain `hp' empirical interpolation.
\newblock {\em International Journal for Numerical Methods in Engineering},
  90:412--428, 2012.

\bibitem{bbFMM}
W.~Fong and E.~Darve.
\newblock The black-box fast multipole method.
\newblock {\em Journal of Computational Physics}, 228(23):8712 -- 8725, 2009.

\bibitem{giebermann}
K.~Giebermann.
\newblock Multilevel approximation of boundary integral operators.
\newblock {\em Computing}, 67(3):183--207, 2001.

\bibitem{coulomb}
Z.~Gimbutas, L.~Greengard, and M.~Minion.
\newblock Coulomb interactions on planar structures: Inverting the square root
  of the laplacian.
\newblock {\em SIAM Journal on Scientific Computing}, 22(6):2093--2108, 2001.

\bibitem{goreinov}
S.A. Goreinov, E.E. Tyrtyshnikov, and N.L. Zamarashkin.
\newblock A theory of pseudoskeleton approximations.
\newblock {\em Linear Algebra and its Applications}, 261(1--3):1 -- 21, 1997.

\bibitem{Greengard}
L.~Greengard and V.~Rokhlin.
\newblock A fast algorithm for particle simulations.
\newblock {\em Journal of Computational Physics}, 135(2):280--292, 1997.

\bibitem{hackbusch3}
W.~Hackbusch.
\newblock A sparse matrix arithmetic based on {H}-matrices. part {I}:
  Introduction to {H}-matrices.
\newblock {\em Computing}, 62(2):89--108, 1999.

\bibitem{hackbusch4}
W.~Hackbusch and S.~B\"orm.
\newblock Data-sparse approximation by adaptive {H}2-matrices.
\newblock {\em Computing}, 69(1):1--35, 2002.

\bibitem{hackbusch}
W.~Hackbusch and S.~B\"orm.
\newblock H2-matrix approximation of integral operators by interpolation.
\newblock {\em Applied Numerical Mathematics}, 43(1--2):129 -- 143, 2002.
\newblock 19th Dundee Biennial Conference on Numerical Analysis.

\bibitem{hackbusch2}
W.~Hackbusch and Z.P. Nowak.
\newblock On the fast matrix multiplication in the boundary element method by
  panel clustering.
\newblock {\em Numerische Mathematik}, 54(4):463--491, 1989.

\bibitem{EIMeff2}
Y.~Maday and O.~Mula.
\newblock A generalized empirical interpolation method: Application of reduced
  basis techniques to data assimilation.
\newblock 4:221--235, 2013.

\bibitem{Maday}
Y.~Maday, N.~C. Nguyen, A.~T. Patera, and S.~Pau.
\newblock A general multipurpose interpolation procedure: the magic points.
\newblock {\em Communications On Pure And Applied Analysis}, 8(1):383--404,
  2008.

\bibitem{genFMM1}
P.G. Martinsson and V.~Rokhlin.
\newblock An accelerated kernel-independent fast multipole method in one
  dimension.
\newblock {\em SIAM J. Sci. Comput}, 26:1389--1404, 2005.

\bibitem{DBLP:M2L}
M.~Messner, B.~Bramas, D.~Coulaud, and E.~Darve.
\newblock Optimized {M2L} {K}ernels for the {C}hebyshev {I}nterpolation based
  {F}ast {M}ultipole {M}ethod.
\newblock {\em CoRR}, abs/1210.7292, 2012.

\bibitem{EIMeff1}
B.~Peherstorfer, D.~Butnaru, K.~Willcox, and H.-J. Bungartz.
\newblock Localized discrete empirical interpolation method.
\newblock {\em SIAM Journal on Scientific Computing}, 36(1):A168--A192, 2014.

\bibitem{tyrtyshnikov}
E.~Tyrtyshnikov.
\newblock Mosaic-skeleton approximations.
\newblock {\em CALCOLO}, 33(1-2):47--57, 1996.

\bibitem{genFMM2}
L.~Ying, G.~Biros, and D.~Zorin.
\newblock A kernel-independent adaptive fast multipole algorithm in two and
  three dimensions.
\newblock {\em J. Comput. Phys.}, 196(2):591--626, 2004.

\bibitem{benchfmm}
R.~Yokota.
\newblock Fast multipole methods benchmark.
\newblock {\em https://sites.google.com/site/rioyokota/research/fmm}.

\end{thebibliography}
\end{document}